\documentclass{amsart}
\input epsf
\newtheorem{thm}{Theorem}[section]
\newtheorem{cor}[thm]{Corollary}
\newtheorem{lem}[thm]{Lemma}
\newtheorem{prop}[thm]{Proposition}
\theoremstyle{definition}
\newtheorem{defn}[thm]{Definition}
\theoremstyle{remark}
\newtheorem{rem}[thm]{Remark}
\numberwithin{equation}{section}

\newcommand{\inv}{^{-1}}
\newcommand{\C}{\mathbb C}
\newcommand{\Z}{\mathbb Z}
\newcommand{\R}{\mathbb R}
\newcommand{\N}{\mathbb N}

\newcommand{\x}{\times}

\newcommand{\Lie}{\operatorname{Lie}}

\newcommand{\Sh}{\operatorname{Sh}}
\newcommand{\hl}{\operatorname{hl}}
\newcommand{\la}{\lambda}
\newcommand{\al}{\alpha}
\newcommand{\si}{\sigma}
\newcommand{\del}{\delta}
\newcommand{\Vand}{\operatorname{Vand}}
\newcommand{\PD}{\operatorname{PD}}

\newcommand{\gl}{\mathfrak{gl}}

\newcommand{\qbinom}[2]{\left [\begin{array}{c} #1\\ #2\end{array}\right]}
\begin{document}

\title[Quantum cohomology and total positivity]
{Quantum cohomology rings of Grassmannians and total positivity}%
\author{Konstanze Rietsch}%
\address{DPMMS, Cambridge University and Newnham College, Cambridge}%
\email{rietsch@dpmms.cam.ac.uk}%

\thanks{This work was begun at the Institute for Advanced Study with the support of
 NSF grant DMS 97-29992. The author is currently supported by EPSRC grant GR/M09506/01}

\date{July 25, 2000, revised September 10, 2001}
\begin{abstract}
We give a proof of a result of D. Peterson's identifying the
quantum cohomology ring of a Grassmannian with the reduced
coordinate ring of a certain subvariety of $GL_n$. The totally
positive part of this subvariety is then constructed and we give
closed formulas for the values of the Schubert basis elements on
the totally positive points. We then use the developed methods to
give a new proof of a formula of Vafa and Intriligator and Bertram
for the structure constants (Gromov--Witten invariants). Finally,
we use the positivity of these Gromov--Witten invariants to prove
certain inequalities for Schur polynomials at roots of unity.
\end{abstract}
\maketitle

\section{Introduction}
The group $U^+$ of unipotent upper--triangular matrices in $GL_n$
have on their coordinate ring a nice basis with positive structure
constants. Namely one has the dual of the classical limit of
Lusztig's geometrically constructed canonical basis of the
quantized enveloping algebra. The existence of this basis has been
closely tied \cite{Lus:TotPos94} to the `totally positive part'
of $U^+$ (the matrices with only nonnegative minors).

In this paper we study certain remarkable subvarieties $\mathcal
V_{d,n}$ of $U^+$ which come up in the stabilizer of a particular
standard principal nilpotent element $e$ as closures of the
$1$--dimensional components under Bruhat decomposition. By a
theorem of Dale Peterson's \cite{Pet:Montreal} the quantum
cohomology rings of Grassmannians may be identified with the
coordinate rings of these varieties. Therefore like $U^+$ itself
these varieties have coordinate rings with `canonical' bases on
them (this time coming from Schubert bases), and with positive
structure constants.

Most of Peterson's results, in particular this one, are
unpublished. But \cite{Kos:QCoh} is at least a reference for the
``Peterson variety'' which is Peterson's approach to encoding all
the quantum cohomology rings of partial flag varieties in one go,
and into which the varieties $\mathcal V_{d,n}$ may be embedded
naturally.

We study the varieties $\mathcal V_{d,n}$ here from the point of
view of explicitly describing their totally positive parts. We
show that the totally positive elements come in a single
one-parameter family, and give closed formulas for the values of
the Schubert basis on it. Furthermore we obtain an upper bound
(coming from this totally positive part) for certain Schur
polynomials evaluated at roots of unity.

Since our main approach is to study the coordinate rings directly
and by elementary means, we also find some new very simple proofs
for results about them implied by the isomorphism with quantum
cohomology. This includes properties of their Schubert bases,
their explicit presentation (and thus Peterson's theorem), and
finally a formula of Vafa and Intriligator and Bertram for the
structure constants.

The paper is organized as follows. After an initial section on the
centralizer of the principal nilpotent element $e$ and its Bruhat
decomposition we focus in on the special subvarieties $\mathcal
V_{d,n}$ mentioned above. In Section~\ref{s:CoordRing} we state
Theorem~\ref{t:CoordRing} on the coordinate ring of $\mathcal
V_{d,n}$ and prove it partially. What makes the variety $\mathcal
V_{d,n}$ so accessible to elementary methods is that its points
are very easy to construct explicitly. This is explained next. In
Section~\ref{s:orthogonality} we prove some orthogonality
formulas for values of Schur polynomials at roots of unity. These
facilitate the final step in the proof of Theorem
\ref{t:CoordRing} and are also the essential ingredient for the
proof of the Bertram--Vafa--Intriligator formula in
Section~\ref{s:VafaInt}. Finally, after a little review of total
positivity in Section~\ref{s:TotPos}, we come to determining the
totally positive part of $\mathcal V_{d,n}$. The totally
nonnegative matrices in $\mathcal V_{d,n}$ have a very beautiful
explicit description, as do the values the Schubert basis
elements take on them. In Section~\ref{s:factorization} we write
down how these totally positive matrices factor into products of
elements of simple root subgroups.

In Section~\ref{s:QCoh} we explain the interpretation of the
coordinate ring of $\mathcal V_{d,n}$ in terms of quantum
cohomology and put into context some of the results encountered in
the earlier sections. We also note that one property of this ring
that we have not been able to derive in an elementary way and
which is not at present known other than following from the
geometric definitions, is the positivity of the structure
constants. (The Vafa--Intriligator formula which computes these
is worse than alternating). We use this positivity property in
the final Section~\ref{s:SchurInequalities} to prove an
inequality for values of Schur polynomials at roots of unity.

{\it Acknowledgements.} I was very fortunate to hear many
inspiring lectures by Dale Peterson
and would like to thank him here. I am also grateful to Shahn
Majid for some useful discussions during the writing of this
paper.

\section{The stabilizer of a principal nilpotent}\label{s:U+e}

\subsection{Preliminaries}\label{s:preliminaries}
We recall some standard facts and notations for $GL_n(\C)$. Let
$B^+, B^-$ be the subgroups of upper--triangular respectively
lower--triangular
matrices in $GL_n(\C)$, and $U^+$ and $U^-$ their unipotent
radicals. Their intersection is the maximal torus $T$ of all
diagonal matrices in $GL_n(\C)$ . Let $e_1,\dotsc, e_{n-1}$ and
$f_1,\dotsc, f_{n-1}$ be the standard Chevalley generators in the
Lie algebras of $U^+$ and $U^-$, respectively. So $e_i$ is the
$n\x n$--matrix with unique non--zero entry $1$ found in position
$(i,i+1)$ and $f_i$ its transpose. Let $X^*(T)$ denote the
character group of $T$, written additively,  and $\triangle\subset
X^*(T)$ the set of roots. Denote by $\Pi=\{\alpha_i\in\triangle |\
i=1,\dotsc , n-1\}$ the usual system of simple roots defined by
$t e_i t\inv=\al_i(t) e_i$ for all $t\in T$. The resulting sets of
positive and negative roots are denoted by $\triangle^+$
respectively $\triangle^-$.

We identify the Weyl group of $GL_n$, the symmetric group,  with
the group $W\subset GL_n(\C)$ of permutation matrices. And we let
$s_1,\dotsc , s_{n-1}$ denote the usual Coxeter generators of $W$.
So $s_j$ corresponds to the adjacent transposition $(j,j+1)$. $W$
acts in the usual way on $X^*(T)$ preserving $\triangle$, and for
any $w\in W$, the length $\ell(w)$ is the number of positive roots
sent to negative by $w$. Let $w_0$ be the order--reversing
permutation, so the longest element in $W$.

\subsection{$(U^+)^e$ and its Bruhat decomposition}\label{s:decomp}

\begin{defn} Let us fix the principal nilpotent element
\begin{equation*}
e=\sum_{i=1}^{n-1} e_i\ \ \in\ \ \gl_n(\C).
\end{equation*}
Let $(U^+)^e:=\{u\in U^+\ |\ u e u\inv=e\}$, the stabilizer of $e$
in $U^+$. This is an abelian subgroup of $U^+$ of dimension
$n-1$. The elements of $(U^+)^e$ are precisely those elements of
$U^+$ of the form
\begin{equation}\label{e:generalu}
u= \begin{pmatrix}
  1   &  x_1         &      x_2          &\dots     &x_{n-2}&x_{n-1} \\
      & 1            &     x_1           &          &       &x_{n-2}\\
      &              & \ddots            &  \ddots   &      &\vdots\\
      &              &                   &\ddots     &x_1    &x_2\\
      &              &                   &           &      1&x_1\\
      &              &                   &          &          &1
 \end{pmatrix}.
\end{equation}

We can thereby explicitly identify the coordinate ring
$\C[(U^+)^e]$ with the polynomial ring $\C[x_1,\dotsc, x_{n-1}]$.
\end{defn}

\subsection{}\label{s:parabolics}
To decompose $(U^+)^e$ by the Bruhat decomposition we need to
look a bit more closely at the Weyl group. Let $K$ be a subset of
$\{1,\dotsc, n-1\}$. Then to $K$ we associate the parabolic
subgroup $W_K=\left < s_k\right >_{ k\in K }$ of $W$ and its
unique longest element $w_K$. We also consider the element
$w^K:=w_0 w_K $, which is the minimal length coset representative
in $ w_0 W_K$. Let $\{\hat k\}$ be the set $\{1,\dotsc,
n-1\}\setminus\{k\}$.

We begin with two lemmas which we learned from Dale Peterson.

\begin{lem}[\cite{Pet:QCoh}]\label{l:PetLem}
The elements $w_K\in W$ defined above are characterized by
\begin{equation*}
\{w_K\in W\ |\ K\subset\{1,\dotsc, n-1\}\}= \{w\in W\ | w\cdot
\Pi\subset (-\Pi)\cup\triangle^+ \}.
\end{equation*}
\end{lem}
\begin{proof}
Suppose that $w=w_K$ is the longest element in $W_K=\left < s_k\
|\ k\in K\right >$. Then $w_K$ has the property that it sends
$\Pi_K=\{\al_k | k\in K\}$ to $-\Pi_K$, while all other simple
roots remain positive. This implies the inclusion $\subseteq$.

Now let $w\in W$ lie in the right--hand side. So $w
\cdot\Pi\subset (-\Pi)\cup\triangle^+$. Then let $K$ be defined by
$K=\{j\ |\ -\al_j\in w\cdot \Pi\}$.  By this definition $w_K\inv
w\cdot\Pi\subset\triangle^+$. Therefore $w=w_K$. Note that
$w_K\inv w\cdot\alpha_i$ is positive also when
$w\cdot\alpha_i\in\Delta^+$ since any positive root sent to
$-\Delta^+$ by $w_K$  must lie in $\sum_{k\in K}\N\ \al_k$. This
cannot be true of $w\cdot\alpha_i$ since $w\inv\cdot \al_k\in
-\Pi$ for any $k\in K$.
\end{proof}

\begin{lem}[\cite{Pet:QCoh}]\label{l:PetLem2}
Bruhat decomposition induces
\begin{equation}\label{e:decomp}
(U^+)^e=\bigsqcup_{ K\subset\{1,\dotsc, n-1\}} (U^+)^e\cap B^- w^K
B^-.
\end{equation}
\end{lem}

\begin{proof} By Bruhat decomposition we can write $u=b_1 w_0 w
b_2$ for some $b_1\in B^-,\ b_2\in U^-$ and $w\in W$. Since
$u\cdot e:=u e u\inv=e$, we have
\begin{equation*}
w b_2\cdot e=  w_0 b_1\inv \cdot e.
\end{equation*}
The right--hand side of this equation is of the form
$\sum_{j=1}^{n-1} m_j\ f_j + x$ for an upper--triangular matrix
$x$ and some $m_j\in\C$. The left--hand side is $w\cdot (e + y)$
for some lower--triangular matrix $y$. It follows from their
equality that $w\cdot\Pi\subset (-\Pi)\cup \triangle^+$. By Lemma
\ref{l:PetLem} we have $w=w_K$ for some $K\subset \{1,\dotsc ,
n-1\}$.
\end{proof}

\begin{defn}
Let $\Delta_j\in\C[(U^+)^e]$ be the top right--hand corner
$(n-j)\x (n-j)$--minor of $u$ in \eqref{e:generalu},
\begin{equation}\label{e:Delta}
 \Delta_j=\det
 \begin{pmatrix} x_{j} &x_{j+1}& & & x_{n-1}\\
 x_{j-1} &x_j &\ddots & &\\
 & \ddots&\ddots &\ddots &\\
 &&\ddots& x_j &x_{j+1}\\
 x_{2j-n+1}&&&x_{j-1}&x_{j}
\end{pmatrix},
\end{equation}
where we set $x_0=1$ and $x_k=0$ for $k<0$.
\end{defn}

We can now give an explicit description of the parts in the
decomposition \eqref{e:decomp} of $(U^+)^e$.
\begin{lem}\label{l:decomp}
For $u\in (U^+)^e$ define $K_u:=\{j\in\{1,\dotsc, n-1\}\ |\
\Delta_j(u)= 0\}$. Then $u\in B^- w^{K_u}B^-$. In particular,
\begin{equation}
(U^+)^e\cap\overline{ B^- w^K B^-}=\{u\in (U^+)^e\ |\
\Delta_j(u)=0\ \text{all $j\in K$} \}.
\end{equation}
\end{lem}
\begin{proof}
By Lemma \ref{l:PetLem2} we have $u\in B^-w^K B^-$ for some
$K\subset \{1,\dotsc , n-1\}$. We need only to show that $K=K_u$.
Let $v_1,\dotsc, v_n$ be the standard basis of $\C^n$ and
$\{v_{i_1}\wedge v_{i_2}\wedge\dotsc\wedge v_{i_d}\ |\ 1\le
i_1<\cdots<i_d\le n\}$ the corresponding basis of
$\bigwedge^d\C^n$. Then $\Delta_j(u)$ is the matrix coefficient
\begin{equation*}
\Delta_j(u)=\left (u\cdot (v_{j+1}\wedge\dotsc\wedge v_{n})\ ,\
v_{1}\wedge\dotsc\wedge v_{n-j}\right )
\end{equation*}
of $u$ expressed in this standard basis. Write $u=b_1 w_0 w_K
b_2$ for some $b_1\in B^-$ and $b_2\in U^-$. Then we have $u\cdot
(v_{j+1}\wedge\dotsc\wedge v_{n})= \pm b_1 w_0 w_{K}\cdot
(v_{j+1}\wedge \dotsc\wedge v_{n})$. Therefore $\Delta_j(u)\ne 0$
precisely if $w_0 w_{K}\cdot (v_{j+1}\wedge \dotsc\wedge
v_{n})=\pm v_{1}\wedge \dotsc\wedge v_{n-j}$, or equivalently if
\begin{equation}\label{e:parabolics}
w_{K}\cdot (v_{j+1}\wedge \dotsc\wedge v_{n})=\pm v_{j+1}\wedge
\dotsc\wedge v_{n}.
\end{equation}
Now consider the maximal parabolic subgroup $W_{\{\hat j\}}$. It
can be characterized as the group of permutation matrices
preserving the subsets $\left \{v_{1},\dotsc, v_{j}\right \}$ and
$\left \{v_{j+1},\dotsc, v_{n}\right \}$ in $\C^n$. Therefore the
condition \eqref{e:parabolics} is equivalent to $w_K\in W_{\{\hat
j\}}$, that is, to $j\notin K$.
\end{proof}

\begin{rem}
The preceding lemmas can be generalized to arbitrary reductive
linear algebraic groups (where the $\Delta_j$'s are replaced with
the corresponding matrix coefficients in the fundamental
representations).
\end{rem}

\section{The variety $\mathcal V_{d,n}$ and its
coordinate ring}\label{s:CoordRing}

We now go over to studying in detail the one--dimensional
components of $(U^+)^e$, corresponding to $K=\{\hat d\}$. These
are very special in that their elements -- and especially their
totally positive elements -- can be constructed explicitly and in
an elementary way. Studying the higher dimensional components
involves much deeper results, and this more general case will be
treated in the forthcoming paper \cite{Rie:PartFlQCoh}.

\begin{defn}[The variety $\mathcal V_{d,n}$.]
Define $\mathcal V_{d,n}$ to be the closed subvariety of
 $GL_n(\C)$ defined as
\begin{equation*}
\mathcal V_{d,n}:=(U^+)^e \cap \overline{B^- w^{\{\hat d\}} B^-},
\end{equation*}
where $w^{\{\hat d\}}$ is the Weyl group element defined in
\ref{s:parabolics}. Let $\C[\mathcal V_{d,n}]$ denote the
coordinate ring of $\mathcal V_{d,n}$ as a reduced affine
algebraic variety.

\end{defn}

\begin{defn}[The ring $\Lambda_{d,n}$]\label{d:Lambda}
For $k=1,\cdots, n$ set $Y_k\in\Z[X_1,\dotsc, X_d]$ to be the the
$k\x k$ determinant
\begin{equation}\label{e:Y_k}
\left |\begin{matrix}
\ X_1 & X_2 &\cdots & &\\
\ 1 &X_1 &\ddots  & &\\
  &\ddots& \ddots&\ddots&\vdots\ \\
  & &\ddots &X_1 &X_2\  \\
  & &        &  1  &X_1\
\end{matrix}\right |.
\end{equation}
In other words, $Y_k=\det\left (X_{i-j+1}\right )_{i,j=1,\dotsc
,k}$, where $X_0:=1$ and $X_k:=0$ if $k\notin\{0,\dotsc, d\}$. We
define $\Lambda_{d,n}$ to be the ring given in terms of generators
and relations by
\begin{equation*}
\Lambda_{d,n}:=\C[X_1,\dotsc, X_d]/\left( Y_{c+1},\dotsc Y_{n-1}
\right ),
\end{equation*}
where $c=n-d$.
\end{defn}

\begin{rem}[The element $q$ and quantum cohomology]\label{r:QCohRing}
Define $q:= (-1)^{d+1} Y_n$ in $\Lambda_{d,n}$. Then we may
identify
\begin{equation*}
\Lambda_{d,n}\cong\C[X_1,\dotsc, X_d,q]/(Y_{c+1},\dotsc,
Y_{n-1},Y_n+(-1)^d q).
\end{equation*}
Therefore the ring $\Lambda_{d,n}$ coincides with
$qH^*(Gr_d(n),\C)$, the quantum cohomology ring  of the
Grassmannian of $d$--planes in $\C^n$, by the presentation found
in \cite{SiTi:QCoh,Witten:VerAlgGras}. See Section~\ref{s:QCoh}
for more on this.
\end{rem}

By a theorem of Dale Peterson's \cite{Pet:Montreal} the quantum
cohomology rings of Grassmannians are identified with coordinate
rings of the varieties $\mathcal V_{d,n}$ (in fact, his statement
is a generalization of this to arbitrary type). Using the known
presentation of these quantum cohomology rings in type $A$,
Peterson's result for that case may be stated as follows.
\begin{thm}[D. Peterson]
\label{t:CoordRing}
 There is an isomorphism
$\Lambda_{d,n}\overset\phi\to \C[\mathcal V_{d,n}]$ which takes
$X_j$ to $x_j$ for $j=1,\dotsc, d$.
\end{thm}

Our first aim is to give a direct proof of this theorem. We will
not make use of any of the properties of quantum cohomology until
Section \ref{s:SchurInequalities}, keeping our treatment
completely elementary. So in particular any facts about
$\Lambda_{d,n}$ which will be required will be proved by hand,
even if they can be deduced from the isomorphism with quantum
cohomology.

\begin{lem}\label{l:variety}
The algebra homomorphism $\C[X_1,\dotsc, X_d]\to \C[\mathcal
V_{d,n}]$ defined by $X_j\mapsto x_j$ is surjective and has
kernel the radical of the ideal generated by the elements
$Y_{c+1},\dotsc, Y_{n-1}$ defined in \ref{d:Lambda} above.
\end{lem}

\begin{proof}
Let $y_1,\dotsc, y_n\in\C[(U^+)^e]$ be defined by $y_k:=\det\left
(x_{i-j+1}\right )_{i,j=1,\dotsc ,k}$, where $x_0:=1$ and $x_j:=0$
if $j\notin\{0,\dotsc, n\}$. It suffices to show that $\mathcal
V_{d,n}$ is the vanishing set of $x_{d+1},\dotsc, x_{n-1},
y_{c+1},\dotsc, y_{n-1}$.

By Lemma \ref{l:decomp}, the variety $\mathcal V_{d,n}$ consists
of the matrices in $(U^+)^e$ for which all the minors $\Delta_j$
with $j\ne d$ vanish (see \eqref{e:Delta}).

Consider the minors
$\Delta_{n-1},\Delta_{n-2},\dotsc,\Delta_{d+1}$. Note that
$\Delta_{n-1}=x_{n-1}$. The vanishing of the minors
$\Delta_{n-1},\dotsc, \Delta_{d+1}$ implies inductively that the
coordinates $x_{d+1},\dotsc, x_{n-1}$ vanish. The converse
implication is immediate.

Let $u$ be the matrix from \eqref{e:generalu}. It is clear that
the inverse matrix to $u$ is given by
\begin{equation*}
u\inv=\begin{pmatrix}
  1   & -y_1    & y_2    & \cdots  & (-1)^{d} y_d    &\cdots   & (-1)^{n-1} y_{n-1}\\
      & 1       & -y_1    &        &                 &  \ddots          &\vdots\\
      &         & \ddots  &\ddots  &                 &                  &  \\
      &         &         &        &                 &                  &(-1)^{d} y_d\\
      &         &         &        &                 &                  &\\
      &         &         &        &                 &               &\vdots      \\
      &         &         &        &                 &          1       & -y_1\\
      &         &         &        &                 &                    &1
\end{pmatrix}.
\end{equation*}
So $\Delta_1(u)=y_{n-1}=(-1)^{n-1}\Delta_{n-1}(u\inv)$. Moreover,
$\Delta_{j}(u)=(-1)^{n-1}\Delta_{n-j}(u\inv)$ for all
$j=1,\dotsc, n-1$ (see \cite{Macd:SFbook} I. (2.9)).

Therefore the conditions $\Delta_{1}(u)=\dotsc=\Delta_{d-1}(u)=0$
are equivalent to $\Delta_{c+1}(u\inv)=\Delta_{c+2}(u\inv)=\dotsc
=\Delta_{n-1}(u\inv)=0$, and these are equivalent to the
conditions $y_{c+1}=\dotsc =y_{n-1}=0$.
\end{proof}

\begin{rem}
It follows from the lemma that any element in $\Lambda_{d,n}$
gives rise to a function on $\mathcal V_{d,n}$, since we now have
a well--defined map $\phi:\Lambda_{d,n}\to\C[\mathcal V_{d,n}]$.
We will make free use of this by letting
\begin{equation*}
X(u):=\phi(X)(u) \qquad \text{for any $X\in\Lambda_{d,n}$ and
$u\in\mathcal V_{d,n}$.}
\end{equation*}
To fully prove Theorem \ref{t:CoordRing}, it remains to show that
$\Lambda^\C_{d,n}$ is reduced. This will be done in Section
\ref{proof:CoordRing}.
\end{rem}

\subsection{Constructing elements in $\mathcal V_{d,n}$}\label{s:variety}
We would now like to construct explicitly elements inside
$\mathcal V_{d,n}$. For $m_1,\dotsc, m_d\in \C$, let
\begin{equation}\label{e:u}
\bar u_n(m_1,\dotsc, m_d):=
\begin{pmatrix}
  1   & m_1       & m_2      &\dots  & m_d    &0      &        &0 \\
      & 1      &   m_1    & \ddots &         &\ddots   &       & \\
      &        & 1       &\ddots  &          &         &\ddots & 0\\
       &       &         & \ddots &          &         &       & m_d\\
      &        &        &    &  \ddots  &   \ddots  & \ddots   &\vdots\\
      &         &       &        &          &   1     & m_1    & m_2\\
      &         &        &      &           &           &1      &m_1 \\
      &        &         &      &            &         &        &1
\end{pmatrix} \in GL_n(\C).
\end{equation}
Consider the roots $z_1,\dotsc, z_d\in\C$ of the polynomial
\begin{equation*}
z^d-m_1 z^{d-1}+\dotsc +(-1)^d m_d=0,
\end{equation*}
counted with multiplicities. So $m_j=E_j(z_1,\dotsc, z_d)$, where
$E_j$ is the $j$--th elementary symmetric polynomial in
$z_1,\dotsc, z_d$. We define
\begin{equation}
u_n(z_1,\dotsc, z_d)=\bar{u}_n(E_1(z_1,\dotsc, z_d),\dotsc,
E_d(z_1,\dotsc , z_d)).
\end{equation}
There is now a nice description for when a matrix $\bar
u_n(m_1,\dotsc, m_d)$ lies in $\mathcal V_{d,n}$ in terms of the
roots $z_1,\dotsc ,z_d$.
\begin{lem}\label{l:explicitVdn}
The matrix $u_n(z_1,\dotsc, z_d)$ for complex numbers $z_1,\dotsc
,z_d$ lies in $\mathcal V_{d,n}$ if and only if either $z_1=\dotsc
=z_d=0$ or $z_1,\dotsc, z_d$ are all distinct and
$z_1^n=\cdots=z_d^n$. All elements of $\mathcal V_{d,n}$ are of
this form. Furthermore, the element $q\in\Lambda_{d,n}$ defined in
Remark \ref{r:QCohRing} evaluates on $u_n(z_1,\dotsc,
z_d)\in\mathcal V_{d,n}$ to $q(u_n(z_1,\dotsc,
z_d))=(-1)^{d+1}z_1^n$.
\end{lem}
The main argument in this lemma could be proved using the
presentation of the quantum cohomology ring via the
Landau--Ginzburg potential, see
\cite{McDSal:QCohBook} Section 8.4, or more originally the work of
Gepner \cite{Gepner:FusRing} for $q=1$. Below is another proof.

\begin{proof}
Let $H_k$ be the $k$--th complete elementary symmetric polynomial.
Then by the Jacobi--Trudi identity, $y_k(u_n(z_1,\dotsc,
z_d))=H_k(z_1,\dotsc, z_d)$. So by Lemma \ref{l:variety} $\mathcal
V_{d,n}$ consists of all $u_n(z_1,\dotsc, z_d)$ for which
\begin{equation*}
H_{c+1}(z_1,\dotsc, z_d)=\dotsc =H_{n-1}(z_1,\dotsc, z_d)=0.
\end{equation*}
Consider the recursion $H_m(x_1,\dotsc x_k)=H_m(x_1,\dotsc,
x_{k-1})+ x_k H_{m-1}(x_1,\dotsc,x_k)$. Therefore
\begin{equation}\label{e:rec}
H_{m-1}(x_1,\dotsc, x_k) =0\ \implies H_{m}(x_1,\dotsc,
x_k)=H_{m}(x_1,\dotsc, x_{k-1}),
\end{equation}
and if the $x_i$ are all nonzero, the opposite implication holds
as well. Now suppose $H_{c+1}(z_1,\dotsc,
z_d)=\cdots=H_{n-1}(z_1,\dotsc, z_d)=0$ for $z_1,\dotsc,
z_d\in\C$. Then the vanishing of $H_{c+1}$ implies
\begin{equation*}
H_{c+2}(z_1,\dotsc,z_d)=H_{c+2}(z_1,\dotsc, z_{d-1}).
\end{equation*}
If $d>2$ this again equals to zero and it follows that
\begin{equation*}
H_{c+3}(z_1,\dotsc, z_d)=H_{c+3}(z_1,\dotsc,
z_{c-1})=H_{c+3}(z_1,\dotsc, z_{d-2}).
\end{equation*}
We repeat this process until we finally get that
\begin{equation*}
H_{n}(z_1,\dotsc, z_d)=H_n(z_1,\dotsc, z_{d-1})=\cdots=
H_n(z_1)=z_1^n.
\end{equation*}
This implies by symmetry that $z_1^n=z_2^n=\dotsc =
z_d^n=H_n(z_1,\dotsc, z_d)$. Suppose one of the $z_i\ne 0$. Then
the same holds for the others. To prove that the $z_i$ must be
distinct, look one step back. We have
$0=H_{n-1}(z_1,\dotsc,z_d)=H_{n-1}(z_k,z_l)$ for any $k\ne
l\in\{1,\dotsc d\}$. But $H_{n-1}(z,z)=nz^{n-1}$ is nonzero for
$z\ne 0$. Hence $z_k\ne z_l$.

Now suppose we are given $d$ distinct $z_1, \dotsc, z_d$ with
$z_1^n=\dotsc = z_d^n$. It remains to check that $H_{c+1},\dotsc
, H_{n-1}$ vanish, or equivalently (since the $z_i$ are nonzero)
that $H_n(z_1,\dots,z_d)=\cdots= H_n(z_2,z_1)=H_n(z_1)=z_1^n$. We
claim that
\begin{equation*}
H_{m+n}(z_1,\dotsc, z_d)=z_1^n H_m(z_1,\dotsc, z_d).
\end{equation*}
This is obvious if $d=1$. For $d>1$ it follows by induction on $d$
using $z_1^n=z_2^n$ and the divided difference identity
\begin{equation*}
H_{n+m}(z_1,\dotsc , z_d)=\frac{H_{n+m+1}(z_1,z_3,\dotsc,
z_d)-H_{n+m+1}(z_2,z_3,\dotsc,z_d)}{z_{1}-z_2}.
\end{equation*}
Now for any $k=1,\dotsc ,d$ we have
\begin{equation*}
\begin{aligned}
H_n(z_1,\dotsc, z_k)&=\frac{H_{n+1}(z_1,z_3,\dotsc,
z_k)-H_{n+1}(z_2,\dotsc,
z_k)}{z_1-z_2}\\
&=z_1^n\frac{H_1(z_1,z_3,\dotsc, z_k)-H_1(z_2,\dotsc , z_k
)}{z_1-z_2}=z_1^n.
\end{aligned}
\end{equation*}
This concludes the proof.
\end{proof}

\section{Orthogonality formulas for Schur polynomials at roots of unity}\label
{s:orthogonality}

We have already encountered the elementary symmetric polynomials
in $d$ variables $E_1,\dotsc, E_d$ and the complete homogeneous
symmetric polynomials $H_1, H_2,\dotsc $. At this point the other
Schur polynomials also enter the picture. Denote by $S_\la$ the
Schur polynomial in $d$ variables corresponding to the partition
$\la$ of at most $d$ parts. We will write such a partition
$\lambda$ as a $d$--tuple of integers $\lambda =(\lambda_1,\dotsc,
\lambda_d)$ with $\lambda_1\ge\lambda_2\ge\cdots \ge \la_d\ge 0$.
As a visual aid partitions are often represented by their Young
diagrams, see \cite{Macd:SFbook}.

There are many definitions of the Schur polynomials. For instance,
the Schur polynomials are the characters of the polynomial
representations of $GL_d$ (symmetric polynomials in the
eigenvalues). We recall also that the Schur polynomials may be
obtained from the elementary and complete homogeneous symmetric
polynomials by the Jacobi--Trudi and dual Jacobi--Trudi
identities: Let $\lambda^t=(\lambda^t_1,\dotsc,\lambda^t_c)$ be
the conjugate partition to $\lambda$ obtained by exchanging rows
and columns in the Young diagram, then
\begin{equation*}
S_{\lambda}=
 \left |\begin{matrix}
 H_{\lambda_1} & H_{\lambda_1+1} &\cdots   &  H_{\lambda_1+d-1}   \\
 H_{\lambda_2-1} &H_{\lambda_2} & \cdots   & H_{\lambda_2+d-2}   \\
               &               &  \ddots &                  \\
   H_{\lambda_d-d+1}            &    \cdots  & \cdots &H_{\lambda_{d}}  \\
\end{matrix}\right |
  =\left |
\begin{matrix}
 E_{\lambda^t_1} & E_{\lambda^t_1+1} &\cdots &   E_{\lambda^t_1+c-1} \\
E_{\lambda^t_2-1} &E_{\lambda^t_2} & \cdots  & E_{\lambda^t_2+c-2}  \\
               &               &  \ddots &                  \\
 E_{\lambda^t_c+c-1}   &    \cdots           & \cdots &E_{\lambda^t_{c}}
\end{matrix}\right |,
\end{equation*}
where $E_k=H_k=0$ for $k<0$. The Schur polynomial $S_\lambda$ is
homogeneous of degree $|\lambda |:=\lambda_1+\cdots + \lambda_d$,
the size of the partition. Also, the Schur polynomials form a
basis of the ring of symmetric polynomials in $d$ variables. For
further background on Schur polynomials we refer to
\cite{Macd:SFbook}.

\begin{defn}[partitions inside a box and $\PD$]
Let $\Sh(d,c)$ denote the set of partitions whose Young diagram
or shape fits into a $d\x c$ box. In other words $\Sh(d,c)$ is the
set of partitions $\lambda=(\lambda_1,\dotsc,\lambda_d)$ of length
at most $d$ such that $\lambda_1\le c$. We use the shorthand
notation $(m^k)$ for the  partition $(m,\dotsc, m,0,\dotsc, 0)$
with $m$ occurring $k$ times. So for example $(c^d)$ is the
longest partition in $\Sh(d,c)$.

For a partition $\lambda\in \Sh(d,c)$ define
\begin{equation*}
\PD(\lambda):=(c-\lambda_d,c-\lambda_{d-1},\dotsc, c-\lambda_1).
\end{equation*}
The notation $\PD$ stands for `Poincar\'e duality' (see Section
\ref{s:QCoh}).
\end{defn}

\begin{defn}[The set $\mathcal I_{d,n}$] We fix the primitive
$n$--th root of unity $\zeta=\exp(\frac{2\pi\imath}n)$. Let
$\zeta^I:=(\zeta^{i_1},\dotsc, \zeta^{i_d})$ be an unordered $d$--tuple of
distinct $n$--th roots of $(-1)^{d+1}$. Then $I=(i_1,\dotsc,i_d)$
may be chosen uniquely such that $-\frac {d-1}2\le i_1<
i_2<\cdots<i_d\le n-\frac {d+1}2$ and the $i_k$'s are all
integers (respectively half--integers) if $d$ is odd (even).
Denote the set of all such $d$--tuples $I$ by $\mathcal I_{d,n}$.

Let $c=n-d$. If $I=(i_1,\dotsc, i_d)\in\mathcal I_{d,n}$, then
denote by $\Hat I=(\hat i_1,\dotsc, \hat i_c)$ the unique
$c$--tuple, $-\frac {d-1}2\le \hat i_1< \hat i_2<\cdots<\hat
i_c\le n-\frac {d+1}2$, such that $\zeta^{i_1}, \dotsc,
\zeta^{i_d},\zeta^{\hat i_1},\dotsc \zeta^{\hat i_c}$ enumerates
all $n$--th roots of $(-1)^{d+1}$.

There is a bijection
\begin{equation}\label{e:bij}
(\ )^t:\mathcal I_{d,n}\to \mathcal I_{c,n}
\end{equation}
which takes $I=(i_1,\dotsc,i_d)$ to $I^t:=(\frac
n2-\hat{i}_c,\dotsc, \frac n2-\hat{ i}_1)$.

We note that this bijection corresponds exactly to transposition
of shapes, $(\ )^t:\Sh(d,c)\to\Sh(c,d)$, after the identification
\begin{equation}\label{e:identification}
\begin{array}{ccc}
\Sh(d,c)&\leftrightarrow &\mathcal I_{d,n}\\
\la=(\la_1,\dotsc,\la_d)&\mapsto &
I_\la=(\frac{d+1}2+\la_d-d,\dotsc,
\frac{d+1}2+\la_1-1).\\
\end{array}
\end{equation}
We have $I_{\la^t}=(I_\la)^t$. Furthermore if $\| I\|:=
\sum_{k=1}^d i_k$ then $\| I_\la \|=|\la |$.
\end{defn}

We will show the following identities for Schur polynomials.
\begin{prop}[Orthogonality formulas]\label{p:orthogonality} Let
 $I,J\in\mathcal I_{d,n}$ and $z_1,\dotsc , z_d, t\in \C$. Then
\begin{equation*}
\begin{array}{ll}
 \text{\rm{(1)}}&\sum_{\la\in\Sh(d,c)} S_{\la}(z_1,\dotsc, z_d)
 S_{\PD(\la)}(\zeta^J)=S_{(c^d)}(\zeta^J)\prod_{k=1}^d
 \prod_{l=1}^c
(1-z_k\zeta^{-\hat j_l}),\\ &\\
 \text{\rm{(2)}}&\sum_{\la\in\Sh(d,c)}
 S_{\la}(t \zeta^I)S_{\PD(\la)}(t \zeta^J)=\delta_{I,J} \frac{n^d\
 S_{(c^d)}(t \zeta^I)}{|\Vand(\zeta^I)|^{2}},\\ &\\
 \text{\rm (3)} &\sum_{\la\in\Sh(d,c)} S_{\la}(
 \zeta^I)\overline{S_{\la}(
 \zeta^{J})}=\delta_{I,J} \frac{n^d\
 }{|\Vand(\zeta^I)|^{2}}.\\ &\\
\end{array}
\end{equation*}
where $|\Vand(\zeta^I)|$ stands for the absolute value of the
Vandermonde determinant $\prod_{k<j}(\zeta^{i_k}-\zeta^{i_j})$ and
the bar in equation (3) stands for complex conjugation.
\end{prop}

\begin{lem}\label{l:SchurVals}
If $\la\in\Sh(d,c)$ and $I\in\mathcal I_{d,n}$ then we have the
following equality of values of Schur polynomials,
\begin{equation}
S_{\la^t}(\zeta^{I^t}) = \frac{S_{\PD(\la)}(\zeta^I)}
{S_{(c^d)}(\zeta^I)}=\overline{S_{\la}(\zeta^{I})}.
\end{equation}
\end{lem}
\begin{proof}
To begin with, recall that $S_{(c^d)}$ is just the character of
the $c$--th power of the determinant representation of $GL_d$.
The right hand side equality follows from
$\overline{S_\la(\zeta^I)}=S_{\la}(\zeta^{-I})$ and the general
formula
\begin{equation*}
S_\la(z_1\inv,\dotsc , z_d\inv)
=\frac{S_{\PD(\la)}(z_1,\dotsc,z_d)}{S_{(c^d)} (z_1,\dotsc ,z_d )}
\end{equation*}
for the character of the dual representation of $GL_d$.

We now prove that $S_{\la^t}(\zeta^{\frac n2-\Hat I})
=S_{\la}(\zeta^{-I})$. Let $1\le k\le d$  and suppose
$\la\in\Sh(d,c)$ is the partition $(1^k):=(1,\dotsc, 1,0,\dotsc,
0)$ with $1$ appearing $k$ times. So $S_\la=E_k$, the $k$--th
elementary symmetric polynomial, and $S_{\la^t}=H_k$, the $k$--th
complete homogeneous symmetric polynomial. Then we have
\begin{equation*}
\begin{aligned}
p(t)&:=\prod_{j=1}^d
 (1+t\zeta^{-i_j})=1+E_1(\zeta^{-I})t+E_2(\zeta^{-I})t^2+\dotsc +
 E_d(\zeta^{-I})t^d,
\\
 r(t)&:=\prod_{j=1}^d \frac{1}{(1-t\zeta^{\frac n2-\hat i_j})}
 =1+H_1(\zeta^{\frac n2-\hat I})t+H_2(\zeta^{\frac n2-\hat I})
 t^2+\dotsc + H_d(\zeta^{\frac n2-\hat I})t^d+ 0+\\
 &\qquad\qquad \qquad\qquad\qquad +0 + H_n(\zeta^{\frac n2-\hat
 I})t^n+
 \text{\rm higher order terms.}&
\end{aligned}
\end{equation*}
And since
\begin{equation*}
r(t)\inv p(t)=\prod_{k=1}^n(1+ t
\zeta^{-\frac{d+1}2+k})=1+(-1)^{d+1} t^n
\end{equation*}
we get $r(t)(1+(-1)^{d+1} t^n)=p(t)$. Comparing the first $d$
coefficients on either side we see that
$E_k(\zeta^{-I})=H_k(\zeta^{\frac n2-\hat I})$. So the required
identity is proved for $\la=(1^k)$, where $k=1,\dotsc, d$. For
all other partitions in $\Sh(d,c)$ the formula follows from this
special case using the Jacobi--Trudi identity.
\end{proof}

\begin{proof}[Proof of Proposition \ref{p:orthogonality}]
Recall the classical identity due to Littlewood,
\begin{equation}\label{e:Littlewood}
\sum_{\lambda\in\Sh(d,c)}S_{\lambda}(z_1,\dotsc , z_d)
S_{\lambda^t}(w_1,\dotsc, w_c)=\prod_{i=1}^d\prod_{j=1}^c (1 + z_i
w_j),
\end{equation}
which is a consequence of the usual orthogonality of Schur
functions (see Macdonald \cite{Macd:SFbook} I.(4.3'), or Fulton
\cite{Fu:YoungTab} Appendix A). Now replace $(w_1,\dotsc, w_c)$ by
$\zeta^{J^t}$ in \eqref{e:Littlewood} and apply Lemma
\ref{l:SchurVals} to get the orthogonality formula (1). Formulas
(2) and (3) follow by furthermore replacing $(z_1,\dotsc , z_d)$
by $\zeta^I$ and checking that
\begin{equation*}
\prod_{k=1}^d
 \prod_{l=1}^c
(1-\zeta^{i_k}\zeta^{-\hat j_l})=
\begin{cases}0 & \text{\ if $I\ne J$}\\
 \frac{n^d}{\left |\Vand(\zeta^I)\right |^2} &\text{\ if $I=J$}.
\end{cases}
\end{equation*}
This is immediate in the case $I\ne J$.  To verify the identity
for $I=J$ compute
\begin{equation*}
\left |\Vand(\zeta^I)\right |^2 \prod_{k=1}^d
 \prod_{l=1}^c
(1-\zeta^{i_k}\zeta^{-\hat i_l}) = \prod_{{m,r\in\{1,\dotsc,d
\},} \ {m\ne r}} (1-\zeta^{i_m}\zeta^{-i_r})\prod_{k=1}^d
 \prod_{l=1}^c
(1-\zeta^{i_k}\zeta^{-\hat i_l}),
\end{equation*}
which is the same as
\begin{equation}\label{e:I=J}
\prod_{k=1}^d \left [\frac{\prod_{r=1}^n(1-z
\zeta^{i_k}\zeta^{r-\frac{d+1}2})}{1-z}\right ]_{z=1}.
\end{equation}
Since this polynomial in $z$ simplifies to
\begin{equation*}
p(z)=\frac{\prod_{l=1}^n(1-z\zeta^{i_k}\zeta^{-\frac{d+1}2+l})}{1-z}
=\frac{1-z^n}{1-z}=1+z+\dotsc + z^{n-1},
\end{equation*}
we have $p(1)=n$ and \eqref{e:I=J} equals to $n^d$.
\end{proof}

\begin{rem}
Passing from `column orthogonality' to `row orthogonality' in
Proposition \ref{p:orthogonality} we also get the formulas
\begin{equation}\label{e:CharOrth}
\frac{1}{n^d} \sum_{I\in \mathcal I_{d,n}} S_{\lambda}(\zeta^I)
 \overline{S_\mu (\zeta^I)}|\Vand(\zeta^I)|^2=\delta_{\lambda,\mu}
\end{equation}
and
\begin{equation}\label{e:PDOrth}
\frac{1}{n^d} \sum_{I\in \mathcal I_{d,n}} S_{\lambda}(\zeta^I)
 S_{\PD(\mu)} (\zeta^I) \frac{|\Vand(\zeta^I)|^2}
 {S_{(c^d)}(\zeta^I)}=\delta_{\lambda,\mu},
\end{equation}
for any $\la,\mu\in\Sh(d,c)$. The first one, \eqref{e:CharOrth},
looks like it should have an explanation on terms of
representation theory of $GL_d$. The second one is related to the
analogue of Poincar\'e duality in the quantum cohomology ring,
see Section \ref{s:qPD}.
\end{rem}

\section{Applications to the ring $\Lambda_{d,n}$}

In this section we will finish the proof of Theorem
\ref{t:CoordRing}, that $\Lambda_{d,n}$ is the coordinate ring of
$\mathcal V_{d,n}$. For a partition $\lambda\in\Sh(d,c)$, let
$s_\lambda\in\Lambda_{d,n}$ be defined by
\begin{equation}\label{e:slambda}
s_{\lambda}=
  \left |
\begin{matrix}
 X_{\lambda^t_1} & X_{\lambda^t_1+1} &\cdots &   X_{\lambda^t_1+c-1} \\
X_{\lambda^t_2-1} &X_{\lambda^t_2} & \cdots  & X_{\lambda^t_2+c-2}  \\
               &               &  \ddots &                  \\
 X_{\lambda^t_c+c-1}   &    \cdots           & \cdots &X_{\lambda^t_{c}}
\end{matrix}\right |.
\end{equation}
Note that this looks like the dual Jacobi--Trudi formula.
Therefore these elements evaluate on a point $u=u_n(z_1,\dotsc,
z_d)$ in $\mathcal V_{d,n}$ to
$s_\lambda(u)=S_{\lambda}(z_1,\dotsc, z_d)$.

\subsection{The quantum Pieri rule }
We will need to know how these $s_\la$ behave under
multiplication by the generators $X_k$ of $\Lambda_{d,n}$. This
is explained by the (dual) quantum Pieri rule stated below. The
quantum Pieri rule was first found by Bertram \cite
{Be:QSchubCalc} who proved it for the ring $qH^*(Gr_d(n))$ using
some geometry. But later Bertram, Ciocan--Fontanine, and Fulton
\cite{BeCioFu:QSchurMult} also gave an algebraic and
combinatorial proof inside the ring $\Lambda_{d,n}$.

\begin{thm}[dual quantum Pieri rule,
\cite{BeCioFu:QSchurMult} Proposition 4.2] Let $\la\in\Sh(d,c)$
and $k\in\{1,\dotsc, d\}$, then the following formula holds in
$\Lambda_{d,n}$,
\begin{equation*}
X_k s_\la=\sum_{\nu} s_\nu + q \sum_{\mu} s_\mu,
\end{equation*}
where first sum is over all $\nu\in\Sh(d,c)$ such that $|\nu|=|\la
|+k$ and $\nu^t_1\ge\la^t_1\ge \nu^t_2\ge\dotsc\ge
\nu^t_c\ge\la^t_c\ge 0$. The second sum is over all $\mu$ such
that $|\mu|=|\la | + k-n$ and $\la^t_1-1\ge\mu^t_1\ge\la^t_2-1\ge
\mu^t_2\ge \dotsc \ge\la^t_d-1\ge\nu^t_l\ge 0$.
\end{thm}
Such a formula for multiplication at $q=1$ also shows up in the
earlier work \cite{Gepner:FusRing} of Gepner coming from physics.
And Gepner notes that this specialized ring is semi--simple. The
final step of our proof of Theorem \ref{t:CoordRing} below uses a
similar insight to Gepner's.

\subsection{Conclusion of the proof of Theorem
\ref{t:CoordRing}}\label{proof:CoordRing}

By Lemma \ref{l:variety} it remains only to prove that the
homomorphism $\phi:\Lambda_{d,n}^\C\to \C[\mathcal V_{d,n}]$ is
injective. In other words, let $p\in \Lambda_{d,n}^\C$, and
suppose the image of $p$ in $\C[\mathcal V_{d,n}]$, again denoted
by $p$, is zero. Then we must show that $p=0$ in
$\Lambda_{d,n}^\C$.

At the outset, $p$ is some polynomial in the $X_k$'s. By repeated
application of the dual quantum Pieri rule we can write $p$ as a
linear combination $p=\sum_{\la\in\Sh(d,c)} p_\la s_\la$, where
$p_\la\in\C[q]$. Then because $p$ vanishes as function on
$\mathcal V_{d,n}$ we get that
\begin{equation}\label{e:expressp}
0=p(u_n(t\zeta^I))=\sum_\la p_\la(t^n)S_\la(t\zeta^I),
\end{equation}
for any element $u_n(t\zeta^I)\in\mathcal V_{d,n}$. Now for any
fixed $t\in\C$ consider the vector $(p_\la(t^n))_{\la}\in
\C^{\Sh(d,c)}$. By \eqref{e:expressp} this vector lies in the
orthogonal complement to the span $\left
<(S_\la(t\zeta^I))_\la\right
>_{I\in\mathcal I_{d,n}}$ (with respect to the standard
inner product on $\C^{\Sh(d,c)}$). But the vectors
$(S_\la(t\zeta^I))_\la$ for varying $I\in \mathcal I_{d,n}$ are
all linearly independent, by the orthogonality formula in
Proposition \ref{p:orthogonality} (2), and therefore span
$\C^{\Sh(d,c)}$. So it follows that the vector
$(p_\la(t^n))_{\la}=0$. Since this holds for all $t\in \C$ we
have that all of the $p_\la=0$ in $\C[q]$. Hence $p=0$ as
required. \qed
\begin{cor}[Schubert basis]\label{c:SchubBasis}
The elements $\{\ q^k s_\lambda\ |\ k\in\Z_{\ge 0}, \
\lambda\in\Sh(d,c)\}$ form a $\Z$--basis of
$\Lambda_{d,n}^\Z=\Z[X_1,\dotsc, X_d]/(Y_{c+1},\dotsc, Y_{n-1})$.
\end{cor}
\begin{proof}
By repeated application of the quantum Pieri rule any element in
$\Lambda_{d,n}^\Z$ is expressible as a $\Z$--linear combination
of elements $q^k s_\lambda$. These are all linearly independent
by exactly the same argument as in the proof of Theorem
\ref{t:CoordRing} above.
\end{proof}

\section{Structure constants of $\Lambda_{d,n}$}\label{s:VafaInt}

In this section we give an elementary proof of a formula for the
structure constants of $\Lambda_{d,n}$. The formula
\eqref{e:VafaInt} below was first noted by Bertram
\cite{Be:QSchubCalc} to hold in the isomorphic ring
$qH^*(Gr_d(n))$ as a consequence of his quantum Giambelli formula
(Theorem \ref{t:SchubSchur}) and the Vafa--Intriligator formula
for Gromov--Witten invariants. The mathematical proof of the
Vafa--Intriligator formula by Siebert and Tian \cite{SiTi:QCoh}
depends on higher--dimensional residues and the geometric
definition of the structure constants. We give a simple proof for
$\Lambda_{d,n}$ below which instead uses Theorem
\ref{t:CoordRing} and the orthogonality formula from Proposition
\ref{p:orthogonality}.

\begin{prop}\label{p:formulas}
Let $\la_1,\dotsc,\la_N\in \Sh(d,c)$, $I\in\mathcal I_{d,n}$ and
$t\in\C^*$. Let $P$ be an arbitrary homogeneous symmetric
polynomial in $d$ variables.
\begin{enumerate}
\item
We have the following formula for $P(t\zeta^I)$ in terms of Schur
polynomials,
\begin{equation*}
\begin{split}
&P(t\zeta^I)=\sum_{\nu\in\Sh(d,c)} m^P_\nu(t) S_\nu(t\zeta^I)\\
&\text{where}\qquad m^P_\nu(t)=\frac{1}{n^d}\sum_{J\in\mathcal
I_{d,n}}P(t\zeta^J)S_\nu(t\inv\zeta^{-J})|\Vand(\zeta^J)|^2.
\end{split}
\end{equation*}
\item
For any homogeneous symmetric polynomial $m$,
\begin{equation*}
\sum_{J\in\mathcal I_{d,n}} m(\zeta^J)=0\quad \text{unless}\quad
\deg (m)\equiv 0\mod n.
\end{equation*}
 In particular, $m^P_\nu=0$ unless $
|\nu|\equiv \deg (P)\mod n$.
\item
Let $p$ be a (non--zero) homogeneous element in $\Lambda_{d,n}$.
Then the coefficients $p_{\la,k}$ of
\begin{equation*}
p=\sum_{\la\in\Sh(d,c),\ k\in\Z_{\ge 0}} p_{\la,k}\ q^k s_\la
\end{equation*}
written out in the basis $\{\ q^k s_\lambda\ |\ k\in\Z_{\ge 0}, \
\lambda\in\Sh(d,c) \}$ of $\Lambda_{d,n}$  are given by the
formula
\begin{equation*}
p_{\la,k}=\begin{cases}
 \frac{1}{n^d}\sum_{J\in\mathcal I_{d,n}} p(u(\zeta^J))
\overline{S_\la(\zeta^{J})} |\Vand(\zeta^J)|^2 &\text{if $\ kn+|\la|=\deg (p)$,}\\
\hskip 3cm 0 & \text{otherwise.}
\end{cases}
\end{equation*}
\end{enumerate}
\end{prop}

\begin{cor}[Bertram--Vafa--Intriligator Formula]\label{c:VafaInt}
The structure constants for multiplication in the basis $\{q^k
s_\lambda\ |\ \lambda\in\Sh(d,c),k\in\mathbb N\}$ of
$\Lambda_{d,n}$,
\begin{equation*}
s_{\la} s_\mu = \sum_{k\in\N,\la_0\in\Sh(d,c)}
<s_{\la},s_{\mu},s_{\nu}>_k q^k s_{\PD(\nu)},
\end{equation*}
are computed by
\begin{equation}\label{e:VafaInt}
<s_{\la}, s_\mu, s_\nu>_k= \frac{1}{n^d}\sum_{J\in\mathcal
I_{d,n}} S_{\la}(\zeta^J)  S_{\mu}(\zeta^J) S_{\nu}(\zeta^J)
\frac{|\Vand(\zeta^J)|^2}{S_{(c^d)}(\zeta^J)}
\end{equation}
whenever $cd+k n=|\la |+ |\mu |+ |\nu |$, and otherwise by
$<s_{\la},s_\mu, s_{\nu}>_k=0$.
\end{cor}

\begin{proof}[Proof of the Corollary] The formula \eqref{e:VafaInt}
is an immediate consequence of
Proposition \ref{p:formulas}(3) and Lemma \ref{l:SchurVals}.
\end{proof}
\begin{proof}[Proof of Proposition \ref{p:formulas}]
For $J\in\mathcal I_{d,n}$ define the symmetric polynomial
$S_J:=\sum_{\nu\in\Sh(d,c)} S_\nu(\zeta^{-J})S_{\nu}$. By Theorem
\ref{p:orthogonality}(3) we have $S_J(\zeta^{I})=\del_{I,J}
\frac{n^d }{|\Vand(\zeta^J)|^2}$. Therefore,
\begin{equation*}
\begin{aligned}
P(t\zeta^I)&=\frac{1}{n^d}\sum_{J\in\mathcal I_{d,n}} P(t \zeta^J)
|\Vand(\zeta^J)|^2 S_J(\zeta^I)=\\
&=\frac{1}{n^d}\sum_{\nu\in\Sh(d,c)}\sum_{J\in\mathcal
I_{d,n}}P(t\zeta^J) |\Vand(\zeta^J)|^2 S_\nu(t\inv \zeta^{-J})
S_\nu(t\zeta^I).
\end{aligned}
\end{equation*}
So (1) is proved.

For (2) notice that $M(t):=\sum_{J\in\mathcal I_{d,n}}
m(t\zeta^J)$ satisfies $M(t)=M(\zeta t)$. Therefore $M$ is a
polynomial in $t^n$, and unless $n$ divides $\deg(m)$ we must
have $M=0$ . Since $(\deg (P)- |\nu|)$ is the degree of $m^P_\nu$,
the rest of (2) also follows.

To prove any identity in $\Lambda_{d,n}$ it suffices, by Theorem
\ref{t:CoordRing}, to check it `on points', that is, evaluated at
the elements of $\mathcal V_{d,n}$. Therefore we can deduce (3)
directly from (1) and (2).
\end{proof}

\begin{rem}
The ring $\Lambda_{d,n}$ specializes to
$\Lambda_{d,n}^{(q=0)}=\C[x_1,\dotsc, x_d]/(y_{c+1},\dotsc, y_n)$,
the usual cohomology ring of the Grassmannian, at $q=0$. The
(images of the) $s_\la$ with $\lambda\in\Sh(d,c)$ now form a
basis of this ring, and the structure constants for this basis
are the well--known Littlewood--Richardson coefficients.
Therefore the formula \eqref{e:VafaInt} for structure constants in
$\Lambda_{d,n}$ gives when $k=0$ a formula for
Littlewood--Richardson coefficients.
\end{rem}

\section{Total Positivity (a quick review)}\label{s:TotPos}
Before proceeding in the next section to determining the totally
positive part of $\mathcal V_{d,n}$ we will give some background
and motivation for doing this.

 Total
positivity, introduced by Schoenberg
in the
1930's  is the study of invertible matrices with exclusively
nonnegative real minors (and nonnegative entries in particular).
We will strictly speaking call such elements of $GL_n$ totally
nonnegative, and, if all minors are in fact strictly positive,
then totally positive. The definition of total positivity is
basis--dependent (or more precisely dependent on a `pinning' of
$GL_n$, see \cite{Lus:TotPos94}), but totally positive matrices
also have nice intrinsic properties. For example, Gantmacher and
Krein showed that all totally positive matrices are regular
semisimple with distinct positive real eigenvalues.

 The totally positive matrices form an open semialgebraic
subset of the real points $GL_n(\R)$ which we denote by
$GL_n(\R_{>0})$ in analogy with the standard notation
$GL_n(\mathbb K)$ for the $\mathbb K$--valued points of $GL_n$
when $\mathbb K$ is a field. We also write $GL_n(\R_{\ge 0})$ for
the totally nonnegative part of $GL_n$.

In the early 1990's Lusztig \cite{Lus:TotPos94} took the theory of
total positivity a big step further and generalized it to all
reductive algebraic groups, at the same time relating it (in the
simply--laced case) to canonical bases and their deep positivity
properties. Let us look in more detail at the special case of the
subgroup $U^+$, to give a taste of what the totally positive part
looks like and how it is related to the canonical basis. We use
the notations introduced in \ref{s:preliminaries}.

Let $U^+(\R)=U^+\cap GL_n(\R)$. We will consider $GL_n(\R)$ and
any subsets as endowed with the usual Hausdorff topology. Inside
$U^+(\R)$ we define $U^+(\R_{\ge 0})$ to be the semigroup of all
totally nonnegative matrices in $U^+$. This is a semigroup since
the totally nonnegative matrices can be characterized as acting
by nonnegative matrices in all exterior powers of the standard
representation of $GL_n$ (with their standard bases).

In the 1950's A. Whitney proved that $U^+(\R_{\ge 0})$ is the
semigroup generated by the elements $\exp(t e_i)$ for
$t\in\R_{\ge 0}$. The following more precise description of
$U^+(\R_{\ge 0})$ was given by Lusztig in \cite{Lus:TotPos94}.
For $w\in W$ let $U^+_w=U^+\cap B^-w B^-$. Then we define
$U^+_w(\R_{>0})$ to be the set of totally nonnegative matrices in
$U^+_w$. For $w$ equal to the simple reflection $s_i$ it is easy
to see that
\begin{eqnarray*}
\R_{>0}&\overset\sim \longrightarrow & \ U^+_{s_i}(\R_{>0}),\\
t\ &\ \mapsto\ &\ x_i(t):=\exp(te_i),
\end{eqnarray*}
and $U^+_{s_i}(\R_{>0})$ is a connected component of
$U^+_{s_i}(\R)\cong \R^*$. Hence the notation $U^+_{s_i}(\R_{>0})$
rather than $U^+_{s_i}(\R_{\ge 0})$.

 Moreover, for general $w$ taking a reduced expression
$w=s_{i_1}\dotsc s_{i_k}$ and multiplying together the
corresponding maps gives a (semi--algebraic) isomorphism
\begin{eqnarray}\label{e:factorization}
\quad (\R_{>0})^k &\overset\sim \longrightarrow &\ U^+_w(\R_{>0}),\\
(t_1,\dotsc , t_k)&\mapsto&\ x_{i_1}(t_1)\cdot \dotsc\cdot
x_{i_k}(t_k),\nonumber
\end{eqnarray}
see \cite{Lus:TotPos94}, Proposition 2.7.

The largest of these real semi--algebraic cells,
$U^+_{w_0}(\R_{>0})$, which corresponds to the longest element in
$W$, is also called the totally positive part of $U^+$ and
denoted by $U^+(\R_{>0})$. It is open in $U^+(\R)$ and its closure
is the set of all totally nonnegative matrices in $U^+$, denoted
by $U^+(\R_{\ge 0})$.

Now let us describe the connection with canonical bases. Consider
the coordinate ring $\C[U^+]$ as the graded dual to the enveloping
algebra $\mathcal U^+$ of the Lie algebra $\Lie(U^+)$. Then
$\C[U^+]$ has a basis $\mathbb B$ given by the dual canonical
basis (obtained from the canonical basis of the dual quantum
enveloping algebra after specializing the quantum parameter to
$1$). In fact it also has a
$\Z$--form which is spanned by this canonical basis. And the
products of basis vectors are nonnegative linear combinations of
basis vectors: the structure constants are nonnegative integers.
This last point is very deep as it follows from Lusztig's
geometric construction of the coproduct in the quantized
universal enveloping algebra \cite{Lus:QuaGrBook}. And this
positivity property enters crucially into the proof of the
following theorem.

\begin{thm}[Lusztig \cite{Lus:TotPos94},
see also \cite{Lus:PartFlag} 3.13]
Let $u\in U^+$. Then
\begin{eqnarray*}
& u\in U^+(\R_{>0})\ \iff\ b(u)\in\R_{>0} \ \text{ for all
$b\in\mathbb B$,}&\quad
\\
\text{and}\quad & u\in U^+(\R_{\ge 0})\ \iff\ b(u)\in\R_{\ge 0} \
\text{ for all $b\in\mathbb B$.}&
\end{eqnarray*}
\end{thm}

We will see in the next section exactly similar results  for the
variety $\mathcal V_{d,n}$, where the dual canonical basis is
replaced by the Schubert basis.

\section{The totally positive part of $\mathcal
V_{d,n}$}\label{s:PosPart}

We begin with some general observations about $\mathcal V_{d,n}$.
Firstly $\mathcal V_{d,n}$ is a $1$--dimensional affine variety
with a $\C^*$--action, given explicitly by
\begin{equation*}
t\cdot \bar u_n(m_1,\dotsc, m_d)=\bar u_n(t m_1, t^2 m_2,\dotsc ,
t^d m_d )
\end{equation*}
(see \eqref{e:u}).
We identify $\Lambda_{d,n} $ with $\C[\mathcal V_{d,n}]$ by the
isomorphism of Theorem \ref{t:CoordRing}. Consider the morphism
therefore defined by $q\in\Lambda_{d,n}$,
\begin{equation*}
q:\mathcal V_{d,n}\to\C.
\end{equation*}
This map is $\C^*$--equivariant for the $\C^*$--action on $\C$
given by $t\cdot z=t^n z$, for $t\in\C^*,\ z\in\C$. And the fiber
over $1$ is precisely $\{u_n(\zeta^I)\ |\ I\in \mathcal
I_{d,n}\}$. So $q$ is an $\binom nd$--fold covering over $\C^*$,
ramified at $0$. Also, two elements $u_n(\zeta^I)$ and
$u_n(\zeta^J)$ in the fiber of $q$ over $1$ lie in the same
$\C^*$--orbit if and only if $\zeta^J=\zeta^k \zeta^I$ for some
$k=1,\dotsc, n$. Therefore $\mathcal V_{d,n}$ consists of
precisely $\frac 1n\binom{n}d$ regular $\C^*$--orbits, and one
trivial one (given by the identity matrix).

Consider now the real form $\Lambda_{d,n}^\R= \R[X_1,\dotsc,
X_d]/(Y_{c+1},\dotsc,Y_{n-1})=\Lambda_{d,n}^\Z\otimes_\Z \R$, and
the corresponding set of real points $\mathcal
V_{d,n}(\R)=\mathcal V_{d,n}\cap GL_n(\R)$ of $\mathcal V_{d,n}$.
If $d':=\left[\frac d2\right]$ is the greatest integer part of
$\frac d2$, then there are precisely $\binom{n}{d'}$ real points
in the fiber of $q$ over $1$. So $\mathcal V_{d,n}(\R)$ is made up
of $\binom{n}{d'}$ regular $\R^*$--orbits and the one trivial
one, and $q:\mathcal V_{d,n}(\R)\to \R$ simply looks as indicated
in the figure below.

$$ \epsfbox{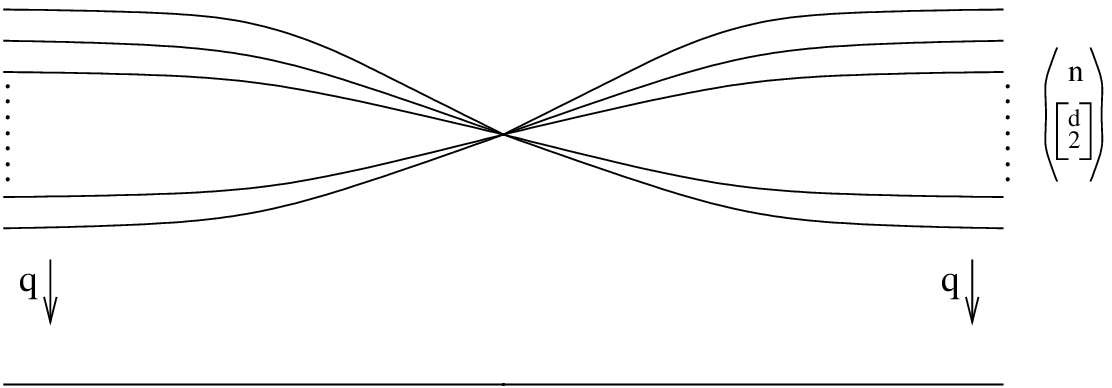} $$

We will identify exactly one of the branches over $\R_{>0}$ as the
totally positive part of $\mathcal V_{d,n}$. Let us give the
precise definitions first.

\begin{defn}
Let $\mathcal V_{d,n}(\R_{\ge 0}):=\mathcal V_{d,n}(\R)\cap
U^+(\R_{\ge 0})$, the set of totally nonnegative matrices in
$\mathcal V_{d,n}$. Recall that by definition
\begin{equation*}
\mathcal V_{d,n}=(U^+)^e\cap\overline{ (B^- w^{\{\hat d\}} B^-)}.
\end{equation*}
Therefore in analogy with the case of $U^+$
described in the previous section we define $\mathcal
V_{d,n}(\R_{>0}):=\mathcal V_{d,n}(\R_{\ge 0})\cap( B^- w^{\{\hat
d\}} B^-)$. The remaining Bruhat decomposition of $\mathcal
V_{d,n}(\R_{\ge 0})$ is very simple since, by Lemma
\ref{l:PetLem2} or the above figure, we have
\begin{equation*}
\mathcal V_{d,n}(\R_{\ge 0})=\mathcal V_{d,n}(\R_{>0})\sqcup
\{id\}.
\end{equation*}
\end{defn}
\begin{defn}
Let $I_0=(-\frac{d-1}{2},-\frac{d-1}{2}+1,\dotsc,
\frac{d-1}{2})\in\mathcal I_{d,n}$. So $\zeta^{I_0}$ is in a
sense the `positive--most' $d$--tuple of roots of $(-1)^{d+1}$.
Then we define a one--parameter family of matrices in $\mathcal
V_{d,n}(\R)$ by
\begin{equation}\label{e:u>0}
 u_{>0}:\R\to \mathcal V_{d,n}(\R), \quad \text{where} \quad
u_{>0}(t):=u_n(t \zeta^{I_0}).
\end{equation}
Note that $u_n(t \zeta^{I_0})$ does indeed lie in $ \mathcal
V_{d,n}(\R)$ since the $d$--tuple $\zeta^{I_0}$ is (up to
permutation) invariant under complex conjugation.
\end{defn}

\subsection{} Let us determine the matrix entries
$E_1(t\zeta^{I_0}),\dotsc , E_d(t\zeta^{I_0})$ of $u_n(t
\zeta^{I_0})$ explicitly. We use a variant of the classical
$q$--binomial theorem of Gauss, with $\zeta$ playing the role of
$q$ (c.f. \cite{Macd:SFbook} I.2 Example 3). For $m,k\in\N$, let
\begin{equation*}
\begin{split}
 [m]_\zeta&:=\frac{\zeta^{\frac{m}2}-\zeta^{-\frac{m}2}}
 {\zeta^{\frac{1}2}-\zeta^{-\frac{1}2}}
 \qquad\qquad\text{and}\\
 \left [\begin{array}{c} m\\ k\end{array}\right ]_\zeta &:=\frac{[m]_\zeta
 [m-1]_\zeta\cdot\dotsc\cdot
 [m-k+1]_\zeta}{[1]_\zeta[2]_\zeta\cdot\dotsc\cdot [k]_\zeta}.
\end{split}
\end{equation*}
Then we have the following version of the binomial formula,
\begin{equation*}
\begin{split}
\left (1+\zeta^{-\frac{d-1}2}t \right )&\left
(1+\zeta^{-\frac{d-1}{2}+1}t\right )\cdot\dotsc\cdot
\left (1+\zeta^{\frac{d-1}2}t\right )=\\
& 1+ {\qbinom d1}_\zeta t + {\qbinom d 2}_\zeta t^2 +\dotsc +
{\qbinom d{d-1}}_\zeta t^{d-1} + t^d.
\end{split}
\end{equation*}
In other words the elementary symmetric polynomials evaluated on
$t\zeta^{I_0}$ give
\begin{equation*}
E_j\left (t\zeta^{-\frac{d-1}2},t\zeta^{-\frac{d-1}2+1},\dotsc,
t\zeta^{\frac{d-1}2}\right )= t^j{\qbinom d j}_\zeta .
\end{equation*}
Since $\zeta=\exp (\frac{2\pi\imath}{n})$ we have that
$\frac{\zeta^{\frac{m}2}-\zeta^{-\frac{m}2}}
{\zeta^{\frac{1}2}-\zeta^{-\frac{1}2}}=\frac{\sin(m\frac \pi
n)}{\sin(\frac \pi n)}$. Therefore
\begin{equation*}
E_j\left (t\zeta^{I_0}\right )=
t^j \frac{\sin(d \frac{\pi}n)\sin((d-1)\frac{\pi}n)
\dotsc\sin((d-j+1)\frac{\pi}n)}{\sin(\frac{\pi}n)\sin(2
\frac{\pi}n)\dotsc \sin(j\frac{\pi}n)}.
\end{equation*}
This gives the $x_j$ entry  of $u_{>0}(t)$.
To determine also the minors of $u_{>0}(t)$ and prove total
positivity we will need a bit more notation.

\begin{defn}
Associate to a partition $\la$ its Young diagram, and let the
boxes in $\la$ be labeled by their $(x,y)$ coordinates in the
usual way, as if they were entries of a matrix. We write
$(i,j)\in\lambda$ if $(i,j)$ are the coordinates of a box in
$\lambda$, so $j\le \lambda_i$ and $i\le \lambda^t_j$. Let $\hl
(i,j)$ denote the {\sl hook length} of the box labeled by
$(i,j)$. That is, $\hl (i,j)=\la_i +\la^t_j-i-j+1$.
\end{defn}

\begin{thm}\label{t:PosPart}
Recall the definitions of $s_\la$ and $q$ (from \eqref{e:slambda}
and Remark \ref{r:QCohRing}).
\begin{enumerate}
\item
For any partition $\la$, the value of $s_\la$ on $u_{>0}(t)$ is
given by the following hook length formula,
\begin{equation*}\label{e:SchurVals}
s_{\lambda}(u_{>0}(t))=t^{|\la |
}\frac{\prod_{(i,j)\in\lambda}\sin((d-i+j)\frac
{\pi}{n})}{\prod_{(i,j)\in\lambda}\sin(\hl(i,j)\frac{\pi}n)}.
\end{equation*}
The value of $q$ is given by $q(u_{>0}(t))=t^n$.
\item
The map $ \R \overset {u_{>0}}\longrightarrow \mathcal
V_{d,n}(\R)$ restricts to give homeomorphisms (for the standard
Hausdorff topology) between
\begin{equation*}
 \R_{\ge 0}\overset\sim \longrightarrow\mathcal V_{d,n}(\R_{\ge
 0})\quad\text{and}\quad  \R_{>0}\overset\sim
 \longrightarrow\mathcal V_{d,n}(\R_{>
 0}).
\end{equation*}
So the $u_{>0}(t)$ for $t\in \R_{\ge 0}$ are precisely all the
totally nonnegative elements in $\mathcal V_{d,n}$.
\item
An element $u\in\mathcal V_{d,n}$ lies in $\mathcal
V_{d,n}(\R_{\ge 0})$ if and only if $s_\la(u)\in\R_{\ge 0}$, for
all $\la\in\Sh(d,c)$. It lies in $\mathcal V_{d,n}(\R_{> 0})$
precisely if $s_\la(u)\in\R_{> 0}$ for all $\lambda\in\Sh(d,c)$.
\end{enumerate}
\end{thm}
\begin{proof}
From the definitions it follows that
$s_\lambda(u_{>0}(t))=S_\lambda(t\zeta^{I_0})=t^{|\lambda
|}S_\lambda(\zeta^{I_0})$. Now the identity in (1) is equivalent
to
\begin{equation*}
S_\lambda(\zeta^{I_0})=\frac{\prod_{(i,j)\in\lambda}\sin((d-i+j)\frac
{\pi}{n})}{\prod_{(i,j)\in\lambda}\sin(\hl(i,j)\frac{\pi}n)}.
\end{equation*}
This equality is again the consequence of a classical formula. It
is a variant of Littlewood's identity for
$S_\la(1,q,q^2,\dotsc,q^k)$, see \cite{Macd:SFbook} Example
I.3.1. The proof goes by writing the Schur polynomial (by Weyl's
character formula) as quotient of two determinants. So
$S_\la(x^{-\frac{d-1}2},x^{-\frac{d-1}2+1},\dotsc,
x^{\frac{d-1}2})=\frac{\Delta^\la_d(x)}{\Delta^{\emptyset}_d(x)}$,
where $\Delta^\la_d(x)$ is the generalized $d\x d$ Vandermonde
determinant
\begin{equation*}
\Delta^\la_d(x)=\det\left( x^{\left (-\frac{d-1}{2}+i-1\right
)(\la_j+d-j)} \right) =\prod_{1\le i<j\le d}\left
(x^\frac{\la_i-\la_j+i-j}{2}-
x^{\frac{\la_i-\la_j+j-i}{2}}\right).
\end{equation*}
Replacing $x$ by $\zeta$ and $\left (x^\frac{\la_i-\la_j+i-j}{2}-
x^{\frac{\la_i-\la_j+j-i}{2}}\right )$ by $2 i\sin\left
((\la_i-\la_j+i-j)\frac{\pi}n\right )$ and cancelling the
numerator against the denominator in $S_\la$ the desired formula
follows. The identity for $q$ is immediate from Lemma
\ref{l:explicitVdn}.

 (2) and (3). To prove that any $u\in U^+$ is totally
nonnegative it suffices to check that all the minors with
connected column sets are nonnegative, by a classical result.
For $u=
u_{>0}(t)=u_n(t\zeta^{I_0})$ these minors are precisely Schur
polynomial values
$S_\la(t\zeta^{-\frac{d-1}2},t\zeta^{-\frac{d-1}2+1},\dotsc,
t\zeta^{\frac{d-1}2})$ where $\lambda$ is a partition which fits
inside a $j\x (n-j$) box for some $j=1,\dotsc n-1$. In other words
$\lambda$ must be a partition with all hook lengths less than $n$.
Therefore in the formula in (3) the denominator is always $\ne 0$
(if $t\ne 0$). From the numerator it follows that
$s_\lambda(u_{>0}(t))\ne 0$ precisely if $\lambda\in\Sh(d,c)$.

Suppose now that $t>0$ and $\lambda\in\Sh(d,c)$. Then for any box
$(i,j)$ in $\la$  both $d-i+j$ as well as the hook length are
positive integers less than $n$. Therefore the values
$\sin((d-i+j)\frac{\pi} n)$ and $\sin(\hl(i,j)\frac{\pi} n)$ in
the formula are all positive and thus these minors are $>0$. For
$t=0$ all $S_\la$ except for $S_{(0,\dotsc , 0)}\equiv 1$ vanish.
This implies (3) and one direction in (2).

It remains to prove that there are no other matrices in $\mathcal
V_{d,n}(\R_{\ge 0})$. Suppose therefore $u=u_n(t\zeta^I)$ is
totally nonnegative, where $t\in\C^*$ and $I\in\mathcal I_{d,n}$.
We have observed that all of the $s_\la$ for $\la\in\Sh(d,c)$
must take nonnegative values on $u$. Since
$u_n(t\zeta^I)\in\mathcal V_{d,n}(\R)$ we know that
$q(u_n(t\zeta^I))=t^n$ must be real. In fact, $q$ can be
expressed as
\begin{equation*}
q=X_d Y_c
\end{equation*}
in $\Lambda_{d,n}$, which follows by expanding the determinant
\eqref{e:Y_k} for $Y_n$ to
\begin{equation*}
Y_n=\sum_{k=1}^{n-1}(-1)^{k+1}X_k Y_{n-k}=(-1)^{d+1} X_d Y_c.
\end{equation*}
Therefore $q(u(t\zeta^I))=t^n$ must be positive real. So
$t=\zeta^k r$ for some $k\in\{1,\dotsc, n\}$ and $r\in\R_{>0}$.
We can replace $t$ by $r$ and $\zeta^I$ by
$\zeta^{J}:=\zeta^{I+(k,\dotsc,k)}$ to get
$u_n(t\zeta^I)=u_n(r\zeta^{J})$. If $\zeta^{J}=\zeta^{I_0}$ then
$u_n(r\zeta^{J})=u_{>0}(r)$, so we may assume that this is not the
case. By the orthogonality formula (2) from Theorem
\ref{p:orthogonality} we then have,
\begin{equation*}
\sum_{\la\in\Sh(d,c)}
 S_{\la}(r \zeta^{J})S_{\PD(\la)}(r \zeta^{I_0})=
 0.
\end{equation*} Since all the $S_{\PD(\la)}(r
\zeta^{I_0})$ are strictly positive, the $S_{\la}(r \zeta^{J})$
cannot all be nonnegative. This finishes the proof.
\end{proof}

\section{Factorization of elements in $\mathcal V_{d,n}(\R_{>0})$}
\label{s:factorization}

The totally positive part $\mathcal V_{d,n}(\R_{>0})$ lies in
particular in $U^+_{w^{\{\hat d\}}}(\R_{>0})$. Hence we can
factorize its elements into products of elements
$x_i(t)=\exp(te_i)$ (with $t\in\R_{>0}$) as described in
\eqref{e:factorization}.

\begin{defn}
For $k\le d$ and $m\le c$ let
\begin{equation*}
B(k,m):=\{(i,j)\ | \ d-k+1\le i\le d ,\ 1\le j\le m\}.
\end{equation*}
Then $B(d,c)$ labels the set of boxes in the Young diagram of
$(c^d)$ and we may think of $B(k,m)$ as parameterizing a subset.
We introduce a partial ordering $\preceq$ on $B(k,m)$ by
\begin{equation*}
(i,j)\preceq (i',j')\quad \iff\quad \text{$i\ge i'$ and $j\ge
j'$}.
\end{equation*}
There is a unique smallest element $(d,m)$ and largest element
$(d-k+1,1)$. When taking products over $B(k,m)$ we will consider
the factors ordered compatibly with this partial order (in
increasing order). So for example the product $\prod_{(i,j)\in
B(d,c)} s_{d-i+j}$ of simple reflections in the Weyl group gives
$w^{\{\hat d\}}$. The different allowed orderings of the factors
correspond to different reduced expressions.

For each pair $(k,m)$ we also have a well--defined monomial
\begin{equation*}
e_{(k,m)}=\prod_{(i,j)\in B(k,m)} e_{d-i+j}
\end{equation*}
in the universal enveloping algebra $\mathcal U^+$ of $\Lie(U^+)$.

Next we define
\begin{equation*}
\begin{array}{rccl}
x:&(\C^*)^{B(d,c)} &\to &U^+\cap B^- w^{\{\hat d\}} B^-\\
&\mathbf a=(a_{(i,j)}) &\mapsto & x(\mathbf a):= \prod_{(i,j)\in
B(d,c)} x_{d-i+j} (a_{(i,j)}).
\end{array}
\end{equation*}
So the restriction of $x$ to $(\R_{>0})^{B(d,c)}$ is the
parameterization map  from \eqref{e:factorization}.

Finally, let  $\Delta_{[i,j]}\in\C[U^+]$ to be the minor with row
set $1,\dotsc, 1+j-i$ and column set $i,\dotsc, j$.
\end{defn}

\begin{lem}\label{l:factorization}
Let $u=x(\mathbf a)$ for $\mathbf a=(a_{(i,j)})\in
(\C^*)^{B(d,c)}$.
\begin{enumerate}
\item
The $m\x m$--minor $\Delta_{[k+1,k+m]}$ of $u$ is given by
\begin{equation*}
\Delta_{[k+1, k+m]}(u)=\prod_{(i,j)\in B(k,m)} a_{(i,j)}.
\end{equation*}
\item
The entries of $\mathbf a$ are computed in terms of the minors
$u$ by
\begin{equation*}\label{e:GenCoeff}
a_{(d-k+1,m)}=\frac{\Delta_{[k+1, k+m]}(u) \Delta_{[k,
k+m-2]}(u)}{\Delta_{[k+1, k+m-1]}(u) \Delta_{[k,k+m-1]}(u)},
\end{equation*}
where $k\in\{1,\dotsc, d\}$ and $m\in\{1,\dotsc ,c\}$.
\end{enumerate}
\end{lem}

The entries $a_{(i,j)}$ can also be worked out using Berenstein
and Zelevinski's generalized Chamber Ansatz \cite{BeZel:TotPos}.
This method for factorizing elements of $U^+_w(\R_{>0})$ is more
general in that it works for any $w$. But it gives the entries in
the factorization of an element $x\in U^+_w(\R_{>0})$ in terms of
minors not of $x$ itself but of another matrix $z=\eta_w\inv(x)$,
which in turn can be expressed using minors of $x$. In the case
at hand it would also be possible to show that their map $\eta_w$
is essentially the identity on these particular elements. This
explains why the formula \eqref{e:GenCoeff} looks like the one
given by a single application of the Chamber Ansatz.

\begin{proof} (1) Let $v_1,\dotsc, v_n$ be the standard basis
of $\C^n$. Then the minor in question is just the matrix
coefficient
\begin{equation}\label{e:MatrixCoeff}
<u\cdot v_{k+1}\wedge\cdots \wedge v_n,\ v_1\wedge\cdots\wedge
v_m\wedge v_{k+m+1}\wedge\cdots\wedge v_n>
\end{equation}
in the fundamental representation $\bigwedge^{n-k} \C^n$ of
$GL_n(\C)$. Now consider the monomial $e(k,m)=\prod_{(i,j)\in
B(k,m)}e_{d-i+j}$ in $\mathcal U^+$. It has the property that
\begin{equation*}
e(k,m)\cdot v_{k+1}\wedge\cdots \wedge v_n=v_1\wedge\cdots\wedge
v_m\wedge v_{k+m+1}\wedge\cdots\wedge v_n.
\end{equation*}
In fact, $e(k,m)$ is the unique (up to scalar) monomial  in the
$e_i's$ satisfying
\begin{equation*}
<e(k,m)\cdot v_{k+1}\wedge\cdots \wedge v_n, v_1\wedge\cdots\wedge
v_m\wedge v_{k+m+1}\wedge\cdots\wedge v_n>\ne 0.
\end{equation*}
Then embed $u=\prod_{(i,j)\in (c^d)}x_{d-i+j}(a_{(i,j)})$ into the
(graded) completion of the enveloping algebra and expand and
multiply the exponentials  to get an infinite series of monomials
in the $e_i$'s. The coefficient of $\prod_{(i,j)\in
B(k,m)}e_{d-i+j}$ in this series is precisely $\prod_{(i,j)\in
B(k,m)} a_{(i,j)}$. Since this is the only monomial that
contributes to the matrix coefficient \eqref{e:MatrixCoeff} the
statement follows. Part (2) follows directly from (1) by
cancellation.
\end{proof}

\begin{prop}\label{p:factorization}
\begin{enumerate}
\item
If $u=u_{>0}(t)$, then
\begin{equation*}
u=\prod_{(i,j)\in B(d,c)}x_{d-i+j}\left (t \
\frac{\sin((i+j-1)\frac \pi n)}{\sin
((d-i+j)\frac{\pi}{n})}\right ),
\end{equation*}
where the order of multiplication is such that the indices of the
$x_{d-i+j}$ spell a reduced expression of $w^{\{\hat d\}}$.
\item
An element $u\in\mathcal V_{d,n}$ lies in $\mathcal
V_{d,n}(\R_{>0})$ precisely if $s_{(m^k)}(u)>0$ for all
rectangular partitions $(m^k)$ with $m\le c$ and $k\le d$.
\end{enumerate}
\end{prop}
\begin{proof}
(1). The minor $\Delta_{[k+1, k+m]}$ of $u_n(t\zeta^I)$ is
computed by $S_{(m^k)}(\zeta^I)$, by the dual Jacobi--Trudi
formula. For $u_{>0}(t)$ this gives
\begin{equation*}
\Delta_{[k+1, k+m]}(u_{>0}(t))=S_{(m^k)}(t\zeta^{I_0})= t^{k
m}\prod_{(i,j)\in (m^k)}\frac{\sin((d-i+j)\frac \pi
n)}{\sin((\hl(i,j))\frac \pi n)}.
\end{equation*}
Therefore $u_{>0}(t)=\prod_{(i,j)\in B(d,c)}x_{d-i+j}(a_{(i,j)})$
where $a_{(i,j)}$ is computed by
\begin{equation*}
a_{(i,j)}=t\ \frac{\sin((i+j-1)\frac \pi n)}{\sin((d-i+j)\frac \pi
n )}
\end{equation*}
after applying Lemma \ref{l:factorization}\eqref{e:GenCoeff} and
cancelling in the numerator and denominator.

(2) is clear since by \eqref{e:factorization} and
\ref{l:factorization} the totally positive  part of $U^+\cap B^-
w^{\{\hat d\}} B^-$ itself is determined by the inequalities
$\Delta_{[m+1,\dotsc, m+k]}>0$. In fact, this criterion is the
earlier result of Berenstein and Zelevinski \cite{BeZel:TotPos},
which we have just reproved in our special case. Now because
$\left . \Delta_{[m+1,\dotsc, m+k]}\right |_{\ \mathcal
V_{d,n}}=s_{(m^k)}$, we are done.
\end{proof}

Note that the reduced set of $c d$ inequalities given in
Proposition \ref{p:factorization}(2) to describe $\mathcal
V_{d,n}(\R_{>0})$ inside variety $\mathcal V_{d,n}(\R)$ is still
surprisingly large, for a $1$--dimensional variety. It is the same
as the set needed to describe the positive part $U_{w^{\{\hat
d\}}}(\R_{>0})$ in the whole $c d$--dimensional variety
$U_{w^{\{\hat d\}}}$.

\section{Quantum cohomology of
Grassmannians}\label{s:QCoh} This section is a brief collecting
together of some basic facts about the (small) quantum cohomology
ring of a Grassmannian. Quantum cohomology originally comes from
the work of Vafa and Witten and ideas from string theory. Since then it
has had a big impact on mathematics, with much work being done to
make the theory rigorous and expand on it. Our main reference is
Bertram's self--contained mathematical treatment for the
Grassmannian \cite{Be:QSchubCalc}. See also \cite{FuPa:QCoh} and
\cite{McDSal:QCohBook} and references therein for a more general
introduction to the mathematical theory of quantum cohomology (in
the algebraic geometric or symplectic setting respectively).

Let $Gr_d(n)$ denote the variety of all $d$--dimensional subspaces
of $\C^n$. As before we fix $c=n-d$. So $c$ is the codimension
while $d$ is the dimension.
\subsection{Schubert basis}\label{s:Schub}
The quantum cohomology of $Gr_{d}(n)$ is in the first instance a
module over the polynomial ring in one variable $\Z[q]$ defined by
\begin{equation}
qH^*(Gr_d(n))=\Z[q]\otimes_\Z H^*(Gr_{d}(n),\Z).
\end{equation}
Therefore it has a $\Z[q]$--basis given by the classical Schubert
basis of $H^*(Gr_{d}(n),\Z)$. We briefly recall its definition
here.

Consider a flag $V_\bullet$ of subspaces $V_\bullet=(V_1\subset
V_2\subset\dotsc\subset V_n=\C^n)$ and a partition
$\lambda\in\Sh(d,c)$. Then the Schubert variety
$X_\la(V_\bullet)\subset Gr_d(n)$ is defined concretely as
\begin{equation*}
X_{\la}(V_{\bullet})= \{W\in Gr_d(n)\ |\ \dim(W\cap V_{c + i -
\la_i})\ge i \}.
\end{equation*}
$X_{\la}(V_{\bullet} )$ is a closed subvariety of $Gr_d(n)$ of
(complex) codimension $|\la |$. Since any flag can be transformed
into any other by the action of $GL_n(\C)$, the homology class of
$X_{\la}(V_\bullet)$ does not depend on the choice of flag
$V_\bullet$. By the classical result of Ehresmann,
the Schubert varieties give rise to a
basis of the homology of $Gr_d(n)$. Let us denote the Poincar\'e
dual basis of cohomology by $\{\si_\la\ |\ \la \in \Sh(d,c) \}$.
So we have $\si_\la\in H^{2 |\la |}(Gr_d(n))$.

\subsection{Poincar\'e duality}\label{s:PD}
Consider the Poincar\'e duality pairing
\begin{equation*}
(\ ,\ ):H^*(Gr_d(n))\x H^*(Gr_d(n))\to\Z
\end{equation*}
defined by $(\si,\tau)=(\si \cup \tau)[Gr_d(n)]$, where
$[Gr_d(n)]\in H_{2 cd}(Gr_{d}(n))$ is the fundamental class.
We recall that in terms of the Schubert basis elements this
pairing takes the form
\begin{equation}\label{e:PD}
(\si_\la,\si_{\mu})=\del_{\la,\PD (\mu)},
\end{equation}
where  as before $\PD(\mu)=(c-\mu_d,c-\mu_{d-1},\dotsc, c-\mu_1)$.

\subsection{Ring structure}\label{s:ring}
As a ring, $qH^*(Gr_d(n))$ does not agree with
$H^*(Gr_d(n))\otimes_\Z \Z[q]$. The new ring structure is a
$q$--deformation of the old with structure constants for the
Schubert basis given by 3--point Gromov--Witten invariants (see
below).
Let us view our Schubert basis elements $\{\si_\la\ |\ \la \in
\Sh(d,c) \}$ as lying inside $qH^*(Gr_d(n))$ by identifying
$\si_\la$ with $\si_\la\otimes 1$.

Let $\la, \mu,\nu\in \Sh(d,c)$ such that $|\la| + |\mu| + |\nu|=
dc + kn$ ( the dimension of the moduli space of degree $k$
holomorphic maps $\C P^1\to Gr_d(n)$ ), and fix three flags
$V_\bullet, V'_\bullet$ and $V''_\bullet$ in general position.
Then $<\si_{\la},\si_{\mu}, \si_{\nu}>_k$ is defined as the
number of degree $k$ holomorphic curves $\C P^1\to Gr_d(n)$
taking $0,1$, and $\infty\in\C P^1$ to points in
$X_{\la}(V_\bullet)$, $X_{\mu}(V'_\bullet)$ and $
X_{\nu}(V''_\bullet)$, respectively. If $|\la| + |\mu| + |\nu|\ne
dc + kn$ then $<\si_{\la},\si_{\mu}, \si_{\nu}>_k:=0$.

Multiplication in $qH^*(Gr_d(n))$ is defined by
\begin{equation*}
\si_\la \cdot \si_\mu=\sum_{\nu, k }
<\si_\la,\si_\mu,\si_{\PD(\nu)}>_k q^k \si_\nu.
\end{equation*}
It is a remarkable fact proved now in more general contexts by
Ruan--Tian \cite{RuTi:QCoh} and also Kontsevich--Manin
\cite{KoMa:QCoh} that the quantum cohomology product is
associative. For the Grassmannian the associativity of the above
defined multiplication is also derived directly in
\cite{Be:QSchubCalc}. With this product $qH^*(Gr_d(n))$ becomes a
graded ring with grading defined by $\deg(\si_\la)= |\la |$ and
$\deg(q)=n$ (or everything doubled to agree with the grading of
the ordinary cohomology). Note that the structure constants of
this ring are nonnegative integers. Also when $q$ is set to zero
only the $k=0$ terms in the product formula remain, and we
recover the definition of the usual cup--product.

\subsection{Presentation and structure constants}\label{s:presentation}

Recall the definition of the ring $\Lambda_{d,n}^\Z$ from
Definition \ref{d:Lambda} (and Corollary \ref{c:SchubBasis}).

\begin{thm}[\cite{SiTi:QCoh}, \cite{Witten:VerAlgGras}]

The ring $qH^*(Gr_d(n))$ is generated by the $d$ Schubert classes
$\si_{(1)},\si_{(1,1)},\dotsc, \si_{(1^d)}$. The assignment
\begin{equation}\label{e:presentation}
\Lambda_{d,n}^\Z\to qH^*(Gr_d(n)):\quad X_k\mapsto
\si_{(1^k)},\quad (-1)^{d+1}Y_n\mapsto q
\end{equation}
defines an isomorphism of rings.
\end{thm}
It is also known that, as in the classical cohomology, the
determinants $Y_1,\dotsc, Y_c$ map to the Schubert classes
$\si_{(1,0,\dotsc,0)},\dotsc,\si_{(c,0,\dotsc, 0)}$. In fact,
Bertram proved that all the classical Giambelli formulas hold
also in the quantum cohomology ring.
\begin{thm}[quantum Giambelli formula, \cite{Be:QSchubCalc}]\label{t:SchubSchur}
For $\la\in\Sh(d,c)$, the element $s_{\la}\in\Lambda^\Z_{d,n}$
defined in \eqref{e:slambda} is mapped to the (quantum) Schubert
class $\si_{\la}$ under the isomorphism \eqref{e:presentation}.
\end{thm}

With this theorem we recover the Bertram--Vafa--Intriligator
formula for the Gromov--Witten invariants from our algebraic
result Corollary \ref{c:VafaInt}.

\begin{thm}[\cite{Be:QSchubCalc}, Section 5]\label{t:VafaInt}
The Gromov--Witten invariant $<\si_\la,\si_\mu,\si_\nu>_k$ where
$|\la |+|\mu |+|\nu |=nk +cd$ is computed by
\begin{equation*}
<\si_{\la}, \si_\mu, \si_\nu>_k= \frac{1}{n^d}\sum_{J\in\mathcal
I_{d,n}} S_{\la}(\zeta^J)  S_{\mu}(\zeta^J) S_{\nu}(\zeta^J)
\frac{|\Vand(\zeta^J)|^2}{S_{(c^d)}(\zeta^J)}
\end{equation*}
\qed
\end{thm}
See also \cite{SiTi:QCoh} for a more general version the
Vafa--Intriligator formula.

\subsection{Poincar\'e Duality for $qH^*(Gr_d(n)$)}\label{s:qPD}

We end this section with the geometric interpretation of the
`row--orthogonality' formula \eqref{e:PDOrth}. One can define an
analogue of the Poincar\'e duality pairing on $qH^*(Gr_d(n))$.
This is a $\Z[q]$--bilinear map
\begin{equation*}
(\ ,\ )_q: qH^*(Gr_d(n))\x qH^*(Gr_d(n))\to \Z[q]
\end{equation*}
which takes $(\si ,\si')$ to the coefficient of $\si_{(c^d)}$ in
$\si \cdot  \si'$. The pairing $(\ ,\ )_q$ specializes to usual
Poincar\'e duality when $q$ goes to  $0$. For two Schubert basis
elements $\si_\lambda$ and $\si_\mu$ this pairing is simply
computed by
\begin{equation}\label{e:qPD}
(\si_\lambda,\si_\mu)_q= <\si_\lambda,\si_\mu,\si_{\PD(c^d)}>_k
q^k
\end{equation}
if $n$ divides $|\la |+ |\mu|-cd$ with quotient $k$. And it is
zero if $|\la |+ |\mu|-cd$ is not divisible by $n$. Now by Theorem
\ref{t:VafaInt}, equation \eqref{e:qPD} becomes
\begin{equation*}
(\si_\lambda,\si_\mu)_q=\frac{1}{n^d}\sum_{J\in\mathcal I_{d,n}}
S_{\la}(\zeta^J)  S_{\mu}(\zeta^J)
\frac{|\Vand(\zeta^J)|^2}{S_{(c^d)}(\zeta^J)}.
\end{equation*}
Thus the equation \eqref{e:PDOrth} corresponds to the following
geometric statement.
\begin{prop}[\cite{Witten:VerAlgGras} Section 3.1,  \cite{Be:QSchubCalc} Lemma
2.5.] Let $\lambda,\mu\in \Sh(d,c)$. Then the pairing $(\ ,\ )_q$
defined above takes the same form
\begin{equation*}(\si_\lambda,\si_\mu)_q=\delta_{\lambda,\PD(\mu)}
\end{equation*}
as the classical Poincar\'e duality pairing. \qed
\end{prop}

\section{An inequality for Schur polynomial values}
\label{s:SchurInequalities}

In this section we give a little application of the positivity of
structure constants in $\Lambda_{d,n}$. The positivity of
structure constants is one property of $\Lambda_{d,n}$ coming
from its isomorphism with $qH^*(Gr_{d}(n),\C)$ for which we know
no elementary explanation. The Vafa--Intriligator formula does not
obviously give something positive. And there is a generalization
of the Littlewood--Richardson rule to quantum cohomology in
\cite{BeCioFu:QSchurMult}, but it is an alternating formula.

\begin{prop}\label{p:SchurInequality}
Let $I_0=(-\frac{d-1}2,\dotsc,\frac{ d-1}2)$ and $I\in \mathcal
I_{d,n}$ any other element, and let $\zeta=\exp(\frac{2\pi i}n)$
as before. Then for any partition $\lambda\in \Sh(d,c)$,
\begin{equation}\label{e:SchurInequality}
\left |\  S_{\lambda}(\zeta^I)\right |\le
S_{\lambda}(\zeta^{I_0}).
\end{equation}
\end{prop}

\begin{rem} An equivalent formulation of the above result is
to say that
\begin{equation*}
\left |\  S_{\lambda}(\zeta^{i_1},\dotsc,\zeta^{i_d})\right |\le
S_{\lambda}(\zeta^{-\frac{d-1}2},\zeta^{-\frac{d-1}2+1},\dotsc,
\zeta^{\frac{d-1}2})
\end{equation*}
for any $d$--tuple of distinct $n$--th roots of unity
$\zeta^{i_1},\dotsc,\zeta^{i_d}$ and partition $\lambda\in\Sh(d,c)$.
\end{rem}

\begin{proof} Consider
$qH^*(Gr_{d}(n),\C)_{q=1}:=qH^*(Gr_{d}(n),\C)\otimes_{\C[q]}
\C[q]/(q-1)$, the specialization of the quantum cohomology ring
at $q=1$. We may identify this ring with  the ring of functions on
the fiber of $\mathcal V_{d,n}$ over $q=1$, that is with
functions on the finite set $\{u_n(\zeta^I)\ |\ I\in \mathcal
I_{d,n}\}$. We have that $qH^*(Gr_{d}(n),\C)_{q=1}$ is a finite
dimensional vector space with basis given by (restriction of) the
Schubert basis $\si_\lambda$.

Another basis of this space is the set consisting of all
$\si_I:=\sum_{\nu\in\Sh(d,c)}
\overline{S_\nu(\zeta^{I})}\si_{\nu}$ for $I\in\mathcal I_{d,n}$.
By Proposition \ref{p:orthogonality} (3) we have
\begin{equation*}
\si_I(u_n(\zeta^{J}))=\del_{I,J} \frac{n^d }{|\Vand(\zeta^I)|^2},
\end{equation*}
so these are up to scalar the characteristic functions of the
points in the set $\{u_n(\zeta^I)\ |\ I\in \mathcal I_{d,n}\}$.

Now consider the multiplication operator $[\si_\la] :\si\mapsto
\si_\la \si$ on $qH^*(Gr_{d}(n),\C)_{q=1}$.
Then clearly
\begin{equation*}
\si_\la\cdot \si_I=S_{\lambda}(\zeta^I)\si_I.
\end{equation*}
So $\si_I$ is an eigenvector of the multiplication operator
$[\si_\la]$ with eigenvalue $S_\lambda(\zeta^I)$. Moreover the
set $\{\si_I\ |\ I\in \mathcal I_{d,n}\}$ is a complete
eigenbasis, in fact a simultaneous eigenbasis for all the
multiplication operators $[\si]$, where $\si\in
qH^*(Gr_{d}(n),\C)_{q=1}$.

Consider the set $Y$ of all $\sigma\in qH^*(Gr_d(n),\R)_{(q=1)}$
such that all the eigenvalues of $[\si]$ are distinct. Then $Y$
is obtained from  the real vector space $qH^*(Gr_d(n),\R)_{(q=1)}$
by removing certain lower dimensional linear subspaces. Therefore
in particular $Y$ is open dense in $qH^*(Gr_d(n),\R)_{(q=1)}$.

Now recall that the multiplication operators $[\si_{\lambda}]$ on
$qH^*(Gr_{d}(n),\C)_{q=1}$ are given in terms of the Schubert
basis by an $\binom{n} {d}\x \binom{n} {d}$--matrix $A_\lambda$
with nonnegative integer entries (Gromov--Witten invariants). We
can approximate the matrix $A_\lambda$ to arbitrary precision
$\epsilon>0$ by some other nonnegative matrix
\begin{equation*}
A_\lambda^\epsilon:=A_\la+\sum_{\mu\in\Sh(d,c)}\epsilon_\mu
A_\mu,\qquad \text{where $0<\epsilon_\mu<\epsilon$, }
\end{equation*}
which has only simple eigenvalues, by the previous paragraph. We
now apply the following version of the Perron--Frobenius
eigenvalue theorem to $A_\lambda^\epsilon$.

 \vskip .2cm
\noindent{\it Perron--Frobenius Theorem.}(\cite{Minc:NonnegMat} I.
Theorem 4.2) Suppose $A$ is an $m\x m$ matrix with all entries in
$\R_{\ge 0}$.  Then $A$ has an eigenvalue $r$ in $\R_{\ge 0}$
with a nonnegative eigenvector and such that all other
eigenvalues of $A$ have absolute value $\le r$. \vskip .2cm

Let $A=A_\lambda^\epsilon$ chosen as above. Since all of its
eigenvalues are simple, the eigenvectors $\si_I$ are unique up to
scalar. Therefore using Theorem \ref{t:PosPart} it follows that
(up to positive scalar) $\sigma_{I_0}$ is the only nonnegative
eigenvector. So by the Perron--Frobenius theorem its eigenvalue,
$S_\lambda(\zeta^{I_0})+\sum_\mu\epsilon_\mu S_\mu(\zeta^{I_0})$
is the maximal eigenvalue. Now choosing $\epsilon$ sufficiently
small we obtain that $S_\lambda(\zeta^{I_0})$ must be a maximal
eigenvalue of $A_\lambda$. This proves the inequality
\eqref{e:SchurInequality}.
\end{proof}

\bibliographystyle{amsplain}
\bibliography{biblio}

\ifx\undefined\bysame
\newcommand{\bysame}{\leavevmode\hbox to3em{\hrulefill}\,}
\fi
\begin{thebibliography}{10}

\bibitem{BeZel:TotPos}
A.~Berenstein and A.~Zelevinsky, {\em Total positivity in {S}chubert
  varieties}, Comment. Math. Helv. {\bf 72} (1997), 128--166.

\bibitem{Be:QSchubCalc}
A.~Bertram, {\em Quantum {S}chubert calculus}, Advances in Mathematics {\bf
  128} (1997), 289--305.

\bibitem{BeCioFu:QSchurMult}
A.~Bertram, I.~Ciocan-Fontanine, and W.~Fulton, {\em Quantum multiplication of
  {S}chur polynomials}, J. Algebra {\bf 219} (1999), no.~2, 728--746.

\bibitem{Fu:YoungTab}
W.~Fulton, , L.M.S. Student Texts, no.~35, Cambridge Univ. Press, 1997.

\bibitem{FuPa:QCoh}
W.~Fulton and R.~Pandharipande, {\em Notes on stable maps and quantum
  cohomology}, Algebraic geometry---Santa Cruz 1995 (Providence, RI), Amer.
  Math. Soc., Providence, RI, 1997, pp.~45--96.

\bibitem{Gepner:FusRing}
D.~Gepner, {\em Fusion rings and geometry}, Comm. Math. Phys. {\bf 141} (1991),
  381--411.

\bibitem{KoMa:QCoh}
M.~Kontsevich and Yu. Manin, {\em Gromov-{W}itten classes, quantum cohomology,
  and enumerative geometry}, Comm. Math. Phys. {\bf 164} (1994), 525--562.

\bibitem{Kos:QCoh}
B.~Kostant, {\em Flag manifold quantum cohomology, the {T}oda lattice, and the
  representation with highest weight $\rho$}, Selecta Math. (N.S.) {\bf 2}
  (1996), 43--91.

\bibitem{Lus:QuaGrBook}
G.~Lusztig, , Progress in Mathematics, Birkhaeuser, Boston, 1993.

\bibitem{Lus:TotPos94}
\bysame, {\em Total positivity in reductive groups}, Lie theory and geometry:
  in honor of Bertram Kostant (Boston) (G.~I. Lehrer, ed.), Progress in
  Mathematics, vol. 123, Birkhaeuser, Boston, 1994, pp.~531--568.

\bibitem{Lus:PartFlag}
\bysame, {\em Total positivity in partial flag manifolds}, Representation
  Theory {\bf 2} (1998), 70--78.

\bibitem{Macd:SFbook}
I.~G. Macdonald, {\em Symmetric functions and {H}all polynomials}, 2nd ed.,
  Oxford Univ. Press, 1995.

\bibitem{McDSal:QCohBook}
Dusa McDuff and Dietmar Salamon, , University Lecture Series, American
  Mathematical Society, Providence, RI, 1994.

\bibitem{Minc:NonnegMat}
Henryk Minc, {\em Nonnegative matrices}, John Wiley \& Sons Inc., New York,
  1988.

\bibitem{Pet:QCoh}
D.~Peterson, {\em Quantum cohomology of ${G}/{P}$}, Lecture Course, spring
  term, M.I.T., 1997.

\bibitem{Pet:Montreal}
\bysame, {\em Quantum cohomology of ${G}/{P}$}, Seminaire de Mathematiques
  Superieures: Representation Theories and Algebraic Geometry, Universite de
  Montreal, Canada, July~28-Aug.~8 1997 (unpublished lecture notes).

\bibitem{Rie:PartFlQCoh}
K.~Rietsch, {\em Quantum cohomology of partial flag varieties and totally
  positive {T}oeplitz matrices}, preprint.

\bibitem{RuTi:QCoh}
Y.~Ruan and G.~Tian, {\em A mathematical theory of quantum cohomology}, J.
  Diff. Geom. {\bf 42} (1995), 259--367.

\bibitem{SiTi:QCoh}
B.~Siebert and G.~Tian, {\em On quantum cohomology rings of {F}ano manifolds
  and a formula of {V}afa and {I}ntriligator}, Asian J. Math. {\bf 1} (1997),
  no.~4, 679--695.

\bibitem{Witten:VerAlgGras}
E.~Witten, {\em The {V}erlinde algebra and cohomology of the {G}rassmannian},
  Geometry, topology \& Physics, Conf. Proc. Lecture Notes {\bf VI} (1995),
  357--422.

\end{thebibliography}

\end{document}